\documentclass[10pt,twocolumn]{IEEEtran}

\usepackage{amsmath,graphicx,amsfonts,amssymb,amsthm,epsfig,mathrsfs,balance}
\usepackage{subfigure,tikz,color,multirow,algorithm,algpseudocode,algorithmicx}
\usepackage{cite,url,framed,bm,dsfont,lipsum}
\usepackage{cuted}
\newtheorem{remark}{Remark}

\newcommand{\differential}{{\rm{d}}}

\begin{document}

\title{\LARGE \bf Stochastic Uncertainty Propagation in Power System Dynamics using Measure-valued Proximal Recursions}

\author{Abhishek~Halder,~\IEEEmembership{Senior Member,~IEEE,} 
        Kenneth~F. Caluya,
        Pegah~Ojaghi,
        and~Xinbo~Geng,~\IEEEmembership{Member,~IEEE}
\thanks{A. Halder and K.F. Caluya are with the Department of Applied Mathematics, University of California, Santa Cruz, CA 95064, USA, {\texttt{\{ahalder,kcaluya\}@ucsc.edu}}.}
\thanks{P. Ojaghi is with the Department of Computer Science and Engineering, University of California, Santa Cruz, CA 95064, USA, {\texttt{pojaghi@ucsc.edu}}.}%
\thanks{X. Geng is with the School of Electrical and Computer Engineering, Cornell University, Ithaca, NY 14850, USA, {\texttt{xg72@cornell.edu}}.} %
\thanks{This research was partially supported by the NSF award 1923278.}}

\maketitle

\begin{abstract}
We present a proximal algorithm that performs a variational recursion on the space of joint probability measures to propagate the stochastic uncertainties in power system dynamics over high dimensional state space. The proposed algorithm takes advantage of the exact nonlinearity structures in the trajectory-level dynamics of the networked power systems, and is nonparametric. Lifting the dynamics to the space of probability measures allows us to design a scalable algorithm that obviates gridding the underlying high dimensional state space which is computationally prohibitive. The proximal recursion implements a generalized infinite dimensional gradient flow, and evolves probability-weighted scattered point clouds. We clarify the theoretical nuances and algorithmic details specific to the power system nonlinearities, and provide illustrative numerical examples.
\end{abstract}


\section{Introduction}\label{sec:intro}

\IEEEPARstart{S}{tochastic} variabilities in power grid have increased significantly in recent years both in the generation side (e.g., due to growing penetration of renewables) as well as in the load side (e.g., due to widespread adoption of plug-in electric vehicles). Several studies \cite{timko1983monte,nwankpa1992stochastic,odun2012structure,ghanavati2016identifying,apostolopoulou2016assessment} have reported that even small stochastic effects can significantly alter the assessment of transient stability, or the performance of automatic generation control. However, the lack of a rigorous yet scalable stochastic computational framework continues to impede \cite{schwalbe2015mathematical} our ability to perform transient analysis involving time-varying joint probability density functions (PDFs) over the states of a large power system network. In this paper, we present a new algorithm to address this computational need.  

\begin{figure}
\centering
\includegraphics[width=\linewidth]{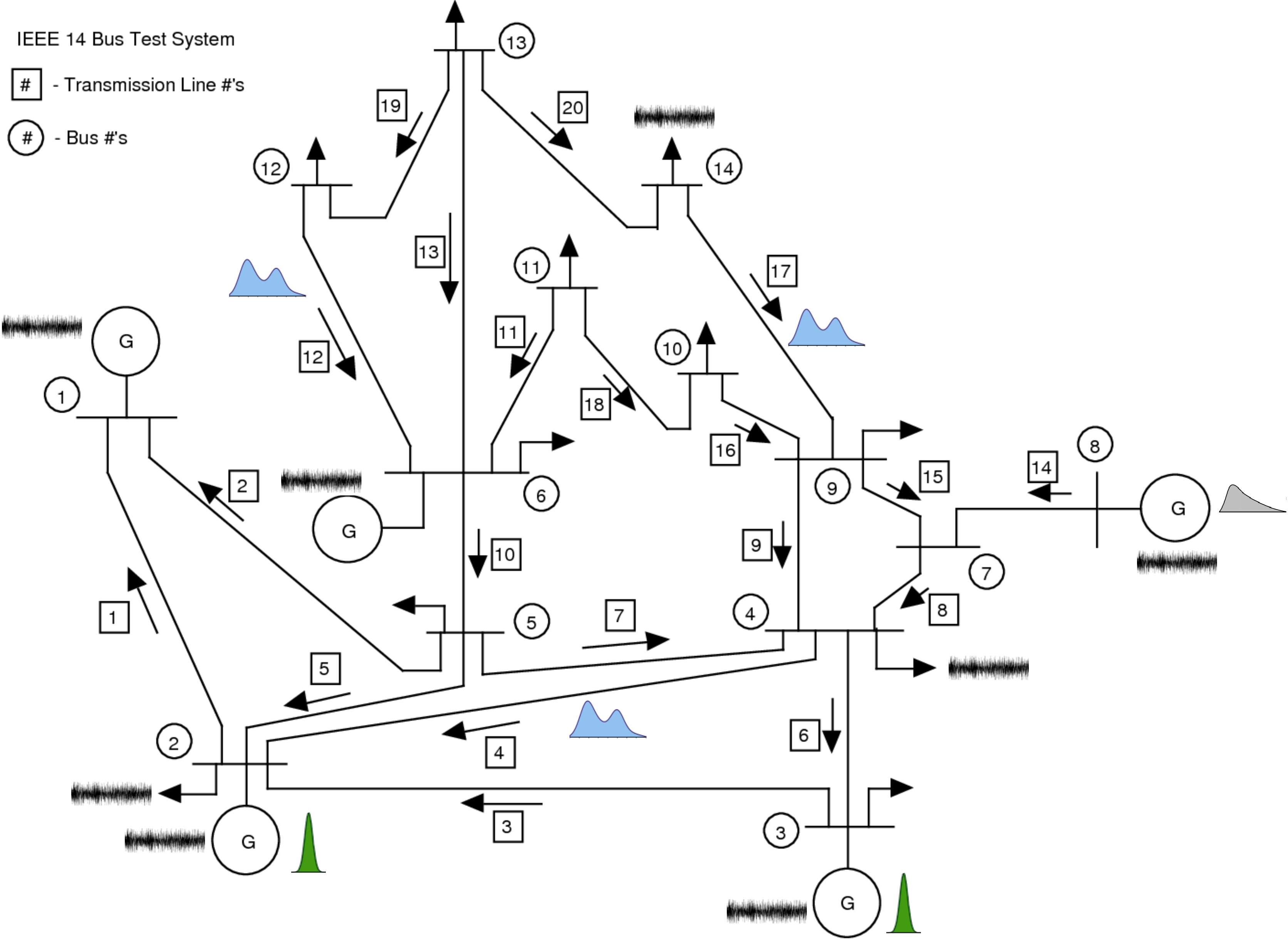}
\caption{\small{A schematic of the IEEE 14 bus test system with stochastic uncertainties. The Uncertainty sources may include stochastic forcing and parametric uncertainties at some generators, random variabilities at some loads, and parametric uncertainties along some transmission lines. For depiction purposes, we indicated the parametric uncertainties as PDFs, and stochastic forcing as intermittent signals.}}
\label{fig:Motivating}
\vspace*{-0.1in}
\end{figure}

\begin{figure}
\centering
\includegraphics[width=0.9\linewidth]{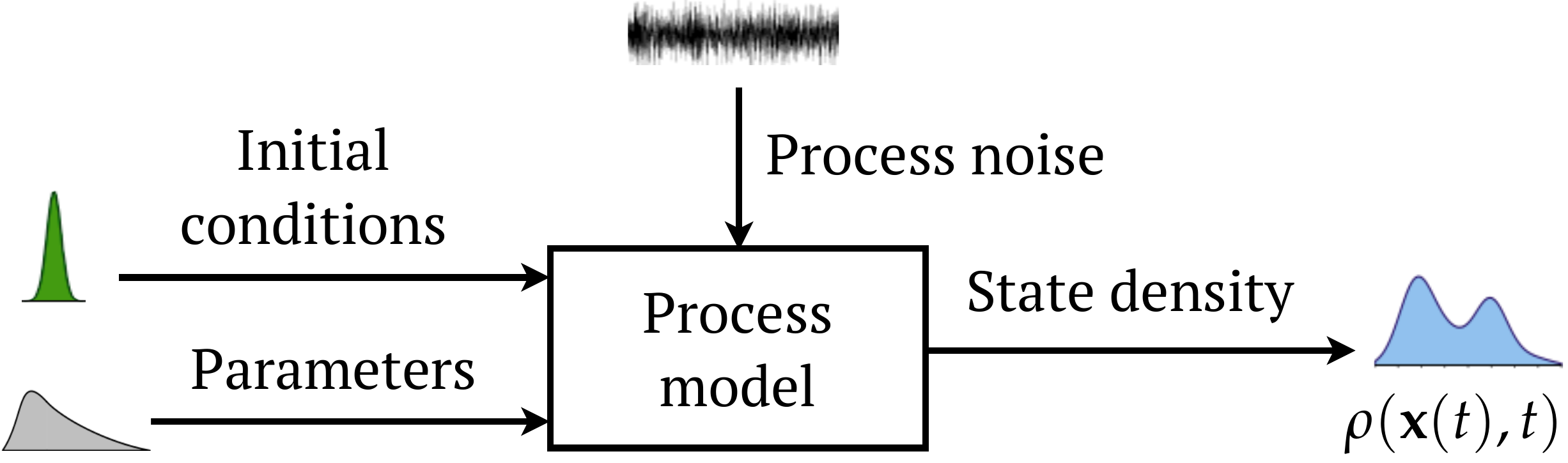}
\caption{\small{Block diagram for joint state PDF propagation.}}
\label{fig:UncertaintyPropagation}
\vspace*{-0.1in}
\end{figure}

Given a networked power system, one can envisage at least three types of uncertainties affecting the dynamics: initial condition uncertainties in the state variables (e.g., rotor phase angles and angular velocities), parametric uncertainties (e.g., inertia and damping coefficients of the generators, reactance associated with different transmission lines), and stochastic forcing (e.g., intermittencies in renewable power generation, load and ambient temperature fluctuations). Fig. \ref{fig:Motivating} depicts a representative scenario. In addition, one could consider uncertainties due to random change in transmission topology resulting from unexpected outage, and uncertainties due to unmodeled dynamics. Given a statistical description of these uncertainties, our approach is to directly solve the \emph{macroscopic flow} of the joint PDFs governing the probabilistic evolution of the state as summarized in Fig. \ref{fig:UncertaintyPropagation}.

\subsubsection{Related works}
Even though the need for quantifying uncertainties in power systems simulations has been long-recognized \cite{allan8791,allan9296}, early studies were limited to statistical reliability assessment. Dynamic simulations with stochastic uncertainties for purposes such as transient stability analysis have been investigated via \emph{Monte Carlo simulations} \cite{timko1983monte,dong2012numerical,perninge2012importance,odun2012structure}.

 As is well known, the Monte Carlo techniques are easy to apply but the computational cost scales exponentially with the number of dimensions, thus making it prohibitive for realistic power systems dynamic simulation. As an alternative, \emph{probabilistic small signal analysis} \cite{nwankpa1992stochastic,rueda2009assessment,huang2013quasi,dhople2013analysis} have appeared in the power systems literature, albeit at the expense of the additional assumption that the random perturbations remain ``small". \emph{Polynomial chaos and related stochastic collocation methods} \cite{hockenberry2004evaluation,xu2019propagating} can do away with the ``small stochastic perturbation" assumption but due to the finite dimensional approximation of the probability space, computational performance degrades if the long-term statistics are desired. Furthermore, to cope with the stochasticity, these techniques require simulating a higher-dimensional system than the dimension of the physical state space, which further limits the scalability for nonlinear simulation. More recently, approximation methods such as \emph{stochastic averaging} \cite{ju2018analytical,ju2018stochastic} based on certain energy function \cite{pai1989energy,chang1995direct,sauerpai1998} have appeared. In \cite{maldonado2018uncertainty}, an algorithm for \emph{propagating first few statistical moments} was proposed. However, this leads to moment closure problems since the dimensions of the time-varying sufficient statistics associated with the corresponding transient joint state PDFs are not known in general. We also note the usage of \emph{stochastic} differential algebraic equations (DAEs) for power system dynamic simulation \cite{wang2011numerical,milano2013systematic}. These studies, however, simulate the sample paths and do not directly propagate the joint state PDFs.

In a different vein, \emph{deterministic bounded uncertainty models} for power flow simulations have been used \cite{dimitrovski2004boundary,hiskens2006sensitivity,chen2012method,althoff2014formal,geng2020reach} for set-valued analysis. These, however, require approximating the underlying nonlinear (deterministic) DAEs appearing in power system dynamic simulation, and thus lead to conservative analysis. For example, the method in \cite{althoff2014formal} requires converting the nonlinear DAEs to linear DAEs in such a way that guarantees set-valued over-approximation of the reachable sets. Likewise, the convex optimization-based bounded uncertainty propagation methods in \cite{choi2017propagating,choi2018propagating} require second order approximation of the power flow state variables as a function of the uncertainties.

\subsubsection{Technical challenges}\label{subsubsec:technicalchallenges}
Several technical challenges need to be overcome to achieve scalable computation enabling the prediction of the joint PDFs over a time horizon of interest. \emph{First}, the trajectory level dynamic models for power systems are inherently nonlinear, which do not preserve Gaussianity, thereby requiring nonparametric prediction of the joint PDFs. \emph{Second}, the joint PDFs for realistic power systems dynamic simulation must evolve over a high dimensional state space, i.e., the joint PDF at any given time has high dimensional support. This necessitates spatial discretization-free algorithms since standard function approximation or interpolation approaches would, in general, be met with ``curse-of-dimensionality" \cite{bellman1957}. The numerical challenges aside, one cannot theoretically guarantee to find a suite of basis functions for the manifold of nonparametric PDFs. \emph{Third}, prediction based on first few statistical moments is challenging since it is not possible to \emph{a priori} guarantee a fixed or even finite dimensional sufficient statistic. For instance, propagating only the mean and covariance could be misleading when the underlying joint PDFs are multi-modal.

\subsubsection{Contributions of this paper}
Our main contribution is to demonstrate that by harnessing recent developments in generalized gradient flows \cite{ambrosio2008gradient} and proximal algorithms \cite{parikh2014proximal}, it is possible to perform \emph{nonparametric} propagation of the joint state PDFs subject to the stochastic nonlinear dynamics of a networked power system, and that the computation can be performed in a gridless and scalable manner. In doing so, the paper also makes theoretical contribution by introducing a change of coordinates that transforms the joint PDF evolution equations in a way that is amenable for the aforesaid generalized proximal recursions.

From a methodological perspective, while proximal algorithms are well-studied \cite{parikh2014proximal,rockafellar1976monotone,combettes2011proximal} in finite dimensional optimization context, the utility of generalized proximal recursions to compute the transient solutions of partial differential equations (PDEs) is relatively less known. This paper highlights this connection by solving the PDE induced by the stochastic power system dynamics.

Our technical approach juxtaposes with the existing power systems literature discussed above, in that we approximate neither the dynamical nonlinearity nor the statistics. Instead of treating the exact nonlinearity as a bane, we exploit the geometry induced by the power system's dynamical nonlinearity over the manifold of time-varying joint state PDFs, thereby enabling a new proximal algorithm to compute the transient joint PDFs.

\subsubsection{Notation and organization}
The set of natural numbers, real numbers, and complex numbers are denoted as $\mathbb{N}$, $\mathbb{R}$, and $\mathbb{C}$, respectively. The symbol $\nabla_{\bm{x}}$ denotes the Euclidean gradient operator with respect to (w.r.t.) the vector $\bm{x}$. Thus, $\nabla_{\bm{x}}\cdot$ stands for the divergence, and $\Delta_{\bm{x}}$ stands for the Laplacian w.r.t. vector $\bm{x}$. We use $\langle\cdot,\cdot\rangle$ to denote the standard Euclidean inner product. The real and imaginary parts of a complex number $z$ are denoted via $\Re(z)$ and $\Im(z)$, respectively. We use the superscript $*$ to denote the complex conjugate, and the superscript $^{\top}$ to denote the matrix transpose. The uniform probability distribution over an interval $[a,b]$ is denoted as ${\rm{Unif}}\left([a,b]\right)$. Likewise, the $n$ dimensional uniform probability distribution over $[a,b]^{n}$ is denoted as ${\rm{Unif}}\left([a,b]^{n}\right)$. The symbol $\propto$ denotes proportionality, $|\cdot|$ denotes the absolute magnitude, $\|\cdot\|_{2}$ denotes the standard Euclidean 2-norm, $\otimes$ denotes the Kronecker product, $\odot$ and $\oslash$ respectively denote the elementwise (Hadamard) product and division, $\det(\cdot)$ stands for the determinant, and the subscript $_{\sharp}$ denotes the pushforward of a PDF via a map. The $n\times n$ identity and zero matrices are denoted as $\bm{I}_{n}$ and $\bm{0}_{n\times n}$, respectively.

The rest of this paper is structured as follows. Section \ref{sec:models} details the power system dynamics models at the microscopic or trajectory level (Sec. \ref{subsec:samplepathdynamics}) as well as at the macroscopic or statistical ensemble level (Sec. \ref{subsec:macroscopic}). The proposed idea of realizing the flow of the joint state PDFs subject to the macroscopic power system dynamics via infinite dimensional proximal recursions, is explained in Section \ref{sec:prox}. Section \ref{sec:ProxAlgo} elucidates the corresponding proximal algorithm that enables the computation of the transient joint state PDFs via weighted point cloud evolution. Numerical simulations illustrating the proposed framework are reported in Section \ref{sec:NumericalSimulations}. Section \ref{sec:conclusion} concludes the paper.


\section{Models}\label{sec:models}

\subsection{Sample Path Dynamics}\label{subsec:samplepathdynamics}

%
%
In this work, we consider the coupled stochastic differential equations (SDEs) associated with the networked-reduced power systems model \cite[Ch. 7]{sauerpai1998}. Specifically, for a power network with $n$ generators, the stochastic dynamics for the $i$-th generator is given by the It\^{o} SDEs
\begin{subequations}
\begin{align}
&{\mathrm{d}}\theta_{i} = \omega_{i}\:{\mathrm{d}}t,\label{RotAngleGen}\\
&m_{i}\:{\mathrm{d}}\omega_{i} =	\left(\!P_{i} - \gamma_{i}\omega_{i} - \displaystyle\sum_{j=1}^{n} k_{ij}\sin\left(\theta_{i}-\theta_{j}-\varphi_{ij}\right)\!\right){\mathrm{d}}t \nonumber\\
&\qquad\qquad\qquad\qquad\qquad\qquad\qquad\qquad\qquad+\sigma_{i}\:{\mathrm{d}}w_{i},
\end{align}
\label{ItoSDEcomponentlevel}
\end{subequations}
where the state variables are the rotor angles $\theta_{i}\in[0,2\pi)$ and the rotor angular velocities $\omega_{i}\in\mathbb{R}$, for $i\in\{1,\hdots,n\}$. The stochastic forcing is modeled through the standard Wiener process $w_{i}(t)$, and the diffusion coefficient $\sigma_{i}>0$ denotes the intensity of stochastic forcing at the $i$th generator.
 
\begin{remark}
We emphasize here that the network-reduced model \eqref{ItoSDEcomponentlevel} is obtained from the so-called structure preserving power network model \cite{odun2012structure} after applying the Kron reduction \cite{dorfler2012kron}, and therefore, has all-to-all connection topology. The derivation of the Kron-reduced parameters is outlined below.
\end{remark}

With the $i$th generator, we associate its inertia $m_{i}>0$ and damping coefficient $\gamma_{i}>0$. The other parameters: the \emph{effective power input} $P_{i}$, the \emph{phase shift} $\varphi_{ij}\in[0,\frac{\pi}{2})$, and the \emph{coupling coefficients} $k_{ij}\geq 0$, depend on the network reduced admittance matrix $\bm{Y}\equiv[Y_{ij}]_{i,j=1}^{n}\in\mathbb{C}^{n\times n}$, $\bm{Y}=\bm{Y}^{\top}$. Specifically,
\begin{subequations}
\begin{align}
 P_{i} &= P_{i}^{\text{mech}} - P_{i}^{\text{load}} - |E_{i}|^{2} \Re\left(Y_{ii}\right) + \Re\left(E_i \cdot I_{i}^{*}\right),\label{defPi}\\
 \varphi_{ij} &= \begin{cases} -\arctan\left(\dfrac{\Re\left(Y_{ij}\right)}{\Im\left(Y_{ij}\right)}\right), & \text{if}\;i\neq j,\\
 0, & \text{otherwise},
 \end{cases}\\
 k_{ij} &= \begin{cases} |E_{i}| |E_{j}| \vert Y_{ij}\vert, & \text{if}\;i\neq j,\\
 0, & \text{otherwise},
 \end{cases}
\end{align}
\label{DefParam}
\end{subequations}
where $P_{i}^{\text{mech}}$ is the mechanical power input, $P_{i}^{\text{load}}$ is the real load, $E_{i}$ is the internal voltage, and $I_{i}\in\mathbb{C}$ is the current for generator $i$.


Suppose that the \emph{unreduced} power network has $n$ generators and $m$ buses. Then the unreduced admittance matrix $\bm{Y}_{\rm{unreduced}}\in \mathbb{C}^{(n+m)\times(n+m)}$ can be partitioned as
\begin{align}
\bm{Y}_{\rm{unreduced}} = \begin{bmatrix} 
\begin{array}{c|c}
\bm{Y}_{\rm{bnd}} & \bm{Y}_{\rm{bnd{\text{-}}int}} \\
\hline\vspace*{-0.12in}\\
\bm{Y}_{\rm{bnd{\text{-}}int}}^{\top} & \bm{Y}_{\rm{int}}
\end{array}	
 \end{bmatrix},
\label{PartionedAdmittanceMatrix}	
\end{align} 
where $\bm{Y}_{\rm{bnd}}\in\mathbb{C}^{n\times n}$, $\bm{Y}_{\rm{int}}\in\mathbb{C}^{m\times m}$, $\bm{Y}_{\rm{bnd{\text{-}}int}}\in\mathbb{C}^{n\times m}$. The matrix \eqref{PartionedAdmittanceMatrix} relates the unreduced current vector $I_{\rm{unreduced}}\in\mathbb{C}^{n+m}$ with the unreduced voltage vector $E_{\rm{unreduced}}\in\mathbb{C}^{n+m}$ via the Kirchhoff equations
\begin{align}
I_{\rm{unreduced}} = \bm{Y}_{\rm{unreduced}} E_{\rm{unreduced}},
\label{DefUnreducedY}	
\end{align}
or equivalently, via its partitioned version associated with the interior and the boundary nodes:
\begin{align}
\begin{array}{c}
\begin{bmatrix}
I_{\rm{bnd}}\\
I_{\rm{int}} 	
\end{bmatrix}	
\end{array} = \bm{Y}_{\rm{unreduced}} \begin{array}{c}
\begin{bmatrix}
E_{\rm{bnd}}\\
E_{\rm{int}} 	
\end{bmatrix}	
\end{array},
\label{PartitionedKirchoffEquation}	
\end{align}
where $I_{\rm{bnd}},E_{\rm{bnd}}\in\mathbb{C}^{n}$ and $I_{\rm{int}},E_{\rm{int}}\in\mathbb{C}^{m}$.

The network reduced admittance matrix $\bm{Y}\in\mathbb{C}^{n\times n}$ in \eqref{DefParam} is the Schur complement of \eqref{PartionedAdmittanceMatrix} w.r.t. the block $\bm{Y}_{\rm{int}}$, i.e.,
\[\bm{Y} = \bm{Y}_{\rm{unreduced}} / \bm{Y}_{\rm{int}} := \bm{Y}_{\rm{bnd}} - \bm{Y}_{\rm{bnd{\text{-}}int}} \bm{Y}_{\rm{int}}^{-1}\bm{Y}_{\rm{bnd{\text{-}}int}}^{\top}.\]
The network reduced current vector $I\in\mathbb{C}^{n}$ in \eqref{defPi} is obtained as
\[I  = \bm{Y}E_{\rm{bnd}} - I_{\rm{bnd}}.\]

One can view (\ref{ItoSDEcomponentlevel}) as the noisy version of the second order nonuniform Kuramoto oscillator model \cite{dorfler2012synchronization,rodrigues2016kuramoto}, given by
\begin{align}
m_{i}\ddot{\theta}_{i} + \gamma_{i}\dot{\theta}_{i} &= P_{i} - \displaystyle\sum_{j=1}^{n} k_{ij}\sin\left(\theta_{i} - \theta_{j} - \varphi_{ij}\right) \nonumber\\
&\qquad\qquad+ \sigma_{i} \times \text{stochastic forcing},
\label{LangevinForm}	
\end{align}
where the stochastic forcing is standard Gaussian white noise. 

We define the \emph{positive diagonal} matrices
\begin{align*}
\bm{M} &:= {\rm{diag}}\left(m_{1},\hdots,m_{n}\right), \nonumber\\
\bm{\Gamma} &:= {\rm{diag}}\left(\gamma_{1},\hdots,\gamma_{n}\right),\nonumber\\
\bm{\Sigma} &:= {\rm{diag}}\left(\sigma_{1},\hdots,\sigma_{n}\right),
\end{align*}
and rewrite (\ref{ItoSDEcomponentlevel}) as a mixed conservative-dissipative SDE in state vector $\bm{x} := (\bm{\theta},\bm{\omega})^{\top} \in \mathbb{T}^{n} \times \mathbb{R}^{n}$ as
\begin{align}
\begin{pmatrix}
{\mathrm{d}}\bm{\theta}\\
{\mathrm{d}}\bm{\omega}	
\end{pmatrix}
\! = \!\begin{pmatrix}
\bm{\omega}\\
-\bm{M}^{-1}\nabla_{\bm{\theta}}V(\bm{\theta}) -\bm{M}^{-1}\bm{\Gamma}\bm{\omega}  	
\end{pmatrix}{\mathrm{d}}t + \!\begin{pmatrix}
 \bm{0}_{n\times n}\\
 \bm{M}^{-1}\bm{\Sigma}	
 \end{pmatrix}{\mathrm{d}}\bm{w},
\label{ItoSDEvectorlevel}
\end{align}
where 
$\bm{w}\in\mathbb{R}^{n}$ is the standard vector Wiener process, $\mathbb{T}^{n}$ denotes the $n$-torus $[0,2\pi)^{n}$, and the potential function $V : \mathbb{T}^{n} \mapsto \mathbb{R}$ is given by
\begin{align}
V(\bm{\theta}) := -\displaystyle\sum_{i=1}^{n} P_{i}\theta_{i} + \!\displaystyle\sum_{\substack{i,j=1\\i<j}}^{n}\!k_{ij}\left(1 - \cos(\theta_{i}-\theta_{j}-\varphi_{ij})\right).
\label{potential}	
\end{align}
The potential (\ref{potential}) has a natural energy function interpretation and can also be motivated by a mechanical mass-spring-damper analogy \cite{dorfler2013synchronization,ishizaki2018}. 

\subsection{Macroscopic Dynamics}\label{subsec:macroscopic}
Given the sample path dynamics \eqref{ItoSDEcomponentlevel} or equivalently \eqref{ItoSDEvectorlevel}, a prescribed initial joint state PDF 
\begin{align}
\rho_{0}(\bm{x})\equiv\rho(t=t_{0},\bm{\theta}(t_{0}),\bm{\omega}(t_{0}))
\label{FPKInitCond}	
\end{align} 
denoting initial condition uncertainties at time $t=t_{0}$, and prescribed parametric uncertainties given by the joint parameter PDF $\rho_{\text{param}}$, the uncertainty propagation problem calls for computing the transient joint state PDFs $\rho(t,\bm{x})\equiv\rho(t,\bm{\theta},\bm{\omega})$ for any desired time $t \geq t_{0}$, which is a nonnegative function supported on the state space $\mathbb{T}^{n} \times \mathbb{R}^{n}$ satisfying $\int \rho = 1$ for all $t\geq t_{0}$.

The corresponding macroscopic dynamics governing the flow of the joint state PDF $\rho(t,\bm{\theta},\bm{\omega})$ is given by a kinetic Fokker-Planck \cite{villani2006hypocoercive} PDE
\begin{align}
\frac{\partial \rho}{\partial t} &= - \langle \bm{\omega},\nabla_{\bm{\theta}}\rho \rangle  + \nabla_{\bm{\omega}} \cdot \left( \rho \left( \bm{M}^{-1}\bm{\Gamma}\bm{\omega} + \bm{M}^{-1}\nabla_{\bm{\theta}}V(\bm{\theta}) \right.\right.\nonumber\\
&\qquad\qquad\qquad \left.\left. +\frac{1}{2} \bm{M}^{-1}\bm{\Sigma}\bm{\Sigma}^{\top}\bm{M}^{-1} \nabla_{\bm{\omega}} \log \rho \right) \right),
\label{KineticFPK}	
\end{align}
subject to the initial condition \eqref{FPKInitCond} and the joint parameter PDF $\rho_{\text{param}}$. Of course, this subsumes special cases such as when either the initial condition or the parameter vector is deterministic. A direct numerical solution of this PDE initial value problem using conventional discretization (e.g., finite difference) or function approximation techniques will not be scalable in general, as explained in Sec. \ref{subsubsec:technicalchallenges}. In the next Section, we discuss how a measure-valued variational recursion proposed in our recent works \cite{caluya2019ACC,caluya2019TAC,halder2020hopfield,caluya2021TAC,caluya2020reflected} can be employed to address this challenge.

We mention here that (\ref{ItoSDEvectorlevel}) has been used in \cite{ju2018analytical,ju2018stochastic} for uncertainty propagation via stochastic averaging approximation where the univariate energy PDF was proposed as a ``proxy" for the entire joint PDF. Most relevant to our approach in the power systems literature is the work in \cite{wang2013fokker}, which indeed voiced the need for computing the transient joint PDFs but only dealt with the single-machine-infinite-bus case -- simplest ($n=1$) instance of (\ref{ItoSDEvectorlevel}). The resulting bivariate Fokker-Planck PDE in \cite{wang2013fokker} was solved via finite element discretization, and revealed rich stochastic dynamics and nontrivial transient stability aspects even in this simple case. However, it is unreasonable to expect that a finite element discretization, or in fact any spatial discretization scheme to solve (\ref{KineticFPK}) for moderately large $n$ in seconds of computational time, thereby limiting our current ability for realistic power systems simulation with stochastic variability. This calls for fundamentally re-thinking what does it mean to solve the PDE (\ref{KineticFPK}) for dynamics (\ref{ItoSDEvectorlevel}).


\section{Measure-valued Proximal Recursion}\label{sec:prox}
\subsection{Generalized Gradient Descent}\label{subsec:GenGradDescent}
Let $\mathcal{P}_{2}\left(\mathbb{T}^{n} \times \mathbb{R}^{n}\right)$ denote the manifold of joint PDFs supported over the state space $\mathbb{T}^{n} \times \mathbb{R}^{n}$, with finite second raw moments. Symbolically,
\begin{align}
\mathcal{P}_{2}\left(\mathbb{T}^{n} \times \mathbb{R}^{n}\right) := &\bigg\{\rho : \mathbb{T}^{n} \times \mathbb{R}^{n} \mapsto \mathbb{R}_{\geq 0} \mid \int \rho = 1,\nonumber\\
&\hspace*{-0.8in}\int\bm{x}^{\top}\bm{x}\:\rho(\bm{x})\:\differential\bm{x} <\infty\; \text{for all}\;\bm{x}\equiv(\bm{\theta},\bm{\omega})^{\top}\in\mathbb{T}^{n} \times \mathbb{R}^{n} \bigg\}.
\label{DefP2}	
\end{align}
We propose to solve the initial value problem for the PDE \eqref{KineticFPK} by viewing its flow $\rho(t,\bm{\theta},\bm{\omega})$ as the gradient descent of some functional $\Phi:\mathcal{P}_{2}\left(\mathbb{T}^{n} \times \mathbb{R}^{n}\right)\mapsto \mathbb{R}_{\geq 0}$ w.r.t. some distance 
\[{\rm{dist}}:\mathcal{P}_{2}\left(\mathbb{T}^{n} \times \mathbb{R}^{n}\right)\times\mathcal{P}_{2}\left(\mathbb{T}^{n} \times \mathbb{R}^{n}\right)\mapsto\mathbb{R}_{\geq 0}.\]
We now explain this idea in detail.

\begin{figure}[tb]
\centering
\includegraphics[width=0.9\linewidth]{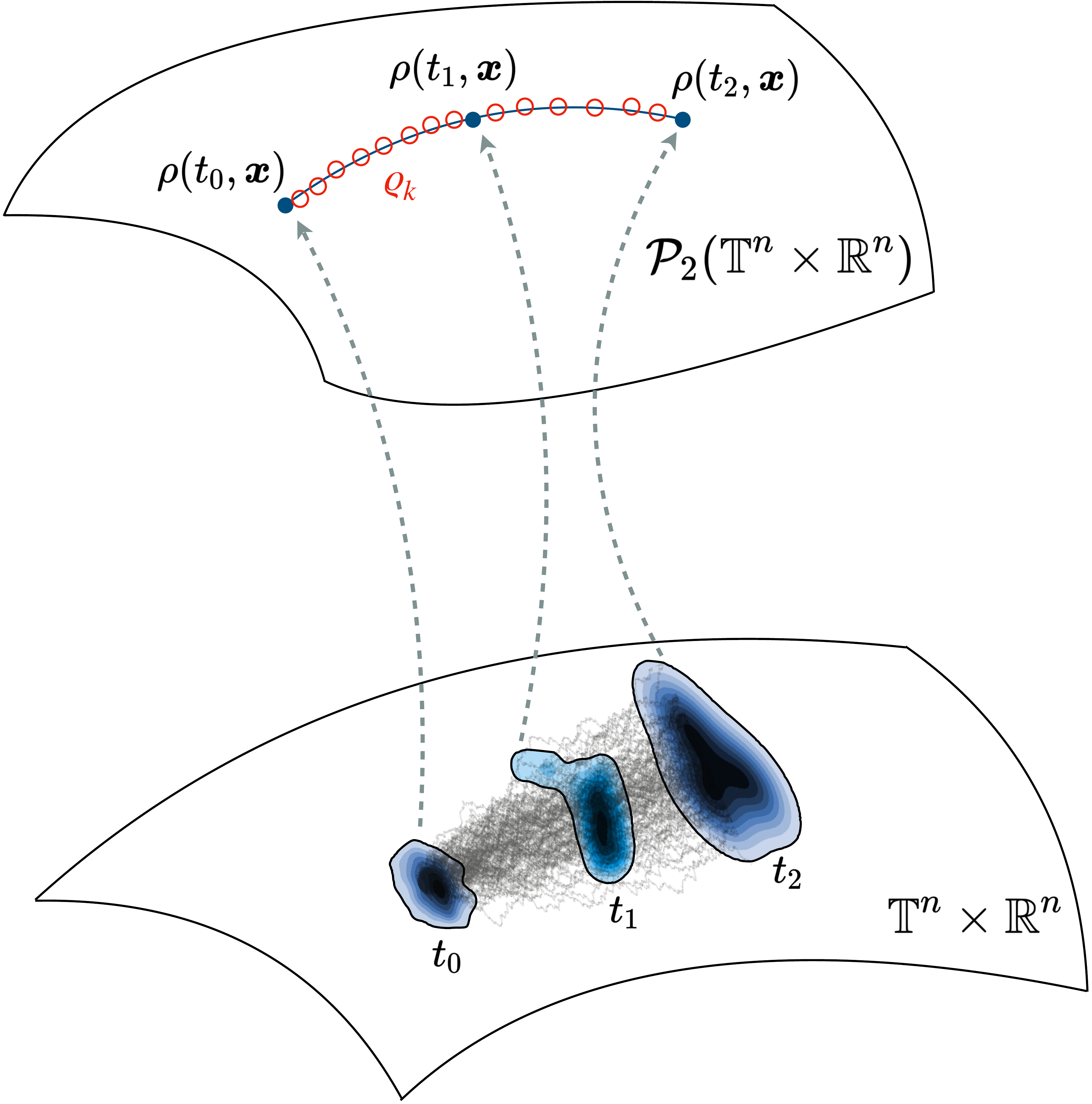}
\caption{\small{A schematic of approximating the flow of the joint PDF trajectory $\rho(t,\bm{x})$ via the sequence $\{\varrho_{k}\}_{k\in\mathbb{N}}$ generated by a variational recursion over $\mathcal{P}_{2}\left(\mathbb{T}^{n}\times\mathbb{R}^{n}\right)$. The curve $\rho(t,\bm{x})$ solves a kinetic Fokker-Planck PDE initial value problem. Points on this curve (shown as three filled circled markers at three specific instances $t_0,t_1,t_2$) are joint PDFs. While the joint PDFs are supported over the finite dimensional base manifold $\mathbb{T}^{n}\times\mathbb{R}^{n}$, the proximal updates $\{\varrho_{k}\}_{k\in\mathbb{N}}$ (shown as circled markers with no face-color) evolve over the infinite dimensional manifold $\mathcal{P}_{2}\left(\mathbb{T}^{n}\times\mathbb{R}^{n}\right)$.}}
\vspace*{-0.1in}
\label{fig:ProxSchematic}
\end{figure}

For $k\in\mathbb{N}$, and for some chosen step size $h>0$, we discretize time as $t_{k} := kh$, and define the infinite dimensional proximal operator of the functional $h\Phi$ w.r.t. the distance ${\rm{dist}}(\cdot,\cdot)$, given by
\begin{eqnarray}
{\text{prox}}^{{\rm{dist}}}_{h\Phi}(\varrho_{k-1}) := \underset{\varrho\in\mathcal{P}_{2}} {{\rm{arg\:inf}}}\: \frac{1}{2} {\rm{dist}}^{2}\left(\varrho,\varrho_{k-1}\right) + h\: \Phi(\varrho).
\label{InfiniteDimProxOpDef}
\end{eqnarray}
Now consider a proximal recursion over the manifold $\mathcal{P}_{2}$ as
\begin{eqnarray}
\varrho_{k} = {\text{prox}}^{{\rm{dist}}}_{h\Phi}(\varrho_{k-1}), \quad k\in\mathbb{N}, \quad \varrho_{0}(\bm{x}) := \rho_{0}(\bm{x}).
\label{ProxRecursionInfiniteDim}	
\end{eqnarray}
Given the PDE \eqref{KineticFPK}, we would like to design the functional pair $\left(\Phi,{\rm{dist}}\right)$ such that the sequence of functions $\{\varrho_{k}\}_{k\in\mathbb{N}}$ generated by the proximal recursion \eqref{ProxRecursionInfiniteDim}, in the small time step limit, converges to the flow $\rho(t=kh,\bm{\theta},\bm{\omega})$ generated by the PDE initial value problem of interest. In particular\footnote{$L^{1}$ denotes the Lebesgue space of absolutely integrable functions.},
\begin{align}
\varrho_{k}(\bm{\theta},\bm{\omega}) \xrightarrow{h\downarrow 0} \rho(t=kh,\bm{\theta},\bm{\omega}) \;\text{in}\;L^{1}\left(\mathbb{T}^{n} \times \mathbb{R}^{n}\right).
\label{ConvergenceGuarantee}	
\end{align}
We remark here that \eqref{DefP2}, \eqref{InfiniteDimProxOpDef}, \eqref{ProxRecursionInfiniteDim}, \eqref{ConvergenceGuarantee} can be written more generally in terms of the joint probability measures instead of PDFs, i.e., even when the underlying measures are not absolutely continuous. Fig. \ref{fig:ProxSchematic} illustrates the idea of approximating the joint PDF trajectory $\rho(t,\bm{x})$ through the sequence $\{\varrho_{k}\}_{k\in\mathbb{N}}$ computed from a variational recursion over $\mathcal{P}_{2}\left(\mathbb{T}^{n}\times\mathbb{R}^{n}\right)$.

We can interpret the proximal operator \eqref{InfiniteDimProxOpDef} as a generalized gradient step for the functional $\Phi$ in the manifold $\mathcal{P}_{2}$. The proximal recursions \eqref{ProxRecursionInfiniteDim} define an infinite dimensional gradient descent of the functional $h\Phi$ over $\mathcal{P}_{2}$ w.r.t. the distance ${\rm{dist}}$. This is reminiscent of the finite dimensional gradient descent, where a gradient flow generated by an ordinary differential equation initial value problem can be recovered as the small time step limit of the sequence of vectors generated by a standard Euclidean proximal recursion; see e.g., \cite[Sec. I]{caluya2019TAC}.

That the flow generated by a Fokker-Planck PDE initial value problem can be recovered from a variational recursion of the form \eqref{ProxRecursionInfiniteDim} was first proposed in \cite{jordan1998variational}, showing that when the drift in the sample path dynamics is a gradient vector field and the diffusion is a scalar multiple of identity matrix, then ${\rm{dist}}(\cdot,\cdot)$ can be taken as the Wasserstein-2 metric arising in the theory of optimal transport \cite{villani2003topics} with $\Phi(\cdot)$ as the free energy functional. In particular, the functional $\Phi$ serves as a Lyapunov functional in the sense $\frac{\differential}{\differential t}\Phi < 0$ along the transient solution of the Fokker-Planck PDE initial value problem. This idea has since been generalized to many other types of PDE initial value problems, see e.g., \cite{ambrosio2008gradient,santambrogio2017euclidean}. 

The algorithmic appeal of the proximal recursion \eqref{ProxRecursionInfiniteDim} is that it opens up the possibility to compute the solution of the PDE initial value problem via recursive convex minimization. A point cloud-based proximal algorithm was proposed in \cite{caluya2019ACC,caluya2019TAC} which was reported to have very fast runtime. The main idea in these references was to \emph{co-evolve} the time-varying state samples (via SDE discretization schemes such as the Euler-Maruyama scheme) as well as the joint PDF values (via \eqref{ProxRecursionInfiniteDim}) evaluated at those samples. The resulting computation is an \emph{online} propagation of the joint PDFs as opposed to the \emph{offline} computation in Monte Carlo or density estimation methods \cite{silverman1998density}. The latter methods only propagate the state samples and then approximate the joint PDFs as post-processing.

Notice that even though the drift in \eqref{ItoSDEvectorlevel} is \emph{not} a gradient vector field, the algorithm in \cite[Sec. V.B]{caluya2019TAC} constructed a pair $\left(\Phi,{\rm{dist}}\right)$ such that \eqref{ProxRecursionInfiniteDim} provably approximates the transient solution of the corresponding kinetic Fokker-Planck PDE with guarantee \eqref{ConvergenceGuarantee}. However, that algorithm cannot be applied to \eqref{KineticFPK} as is. The reasons are explained next.

\subsection{Statistical Mechanics Perspective}
A new difficulty for our SDE \eqref{ItoSDEvectorlevel} is that we have \emph{anisotropic} degenerate diffusion, i.e., the strengths of the noise acting in the last $n$ components of \eqref{ItoSDEvectorlevel} are nonuniform since $\bm{M}^{-1}\bm{\Sigma}$ is not identity. This complicates the matter because the construction of the functional $\Phi$ in \eqref{ProxRecursionInfiniteDim} is usually motivated via free energy considerations utilizing the structure of the \emph{stationary} PDF $\rho_{\infty}(\bm{\theta},\bm{\omega})$ for \eqref{KineticFPK}. The $\rho_{\infty}$ is, in turn, guaranteed to be a unique \emph{Boltzmann distribution} of the form\footnote{here $Z$ is a normalizing constant known as the ``partition function".}
\begin{subequations}
\begin{align}
\rho_{\infty}\left(\bm{\theta},\bm{\omega}\right) &= \dfrac{1}{Z}\exp\left(-\beta H\right), \quad\text{for some}\;\beta > 0,\\
H(\bm{\theta},\bm{\omega}) &:= V(\bm{\theta}) + \frac{1}{2}\langle\bm{\omega},\bm{M\omega}\rangle,	
\end{align}	
\label{BoltzmannPDF}
\end{subequations}
\emph{if and only if} the so-called \emph{Einstein relation} \cite{hernandez1989equilibrium,chen2015fast} holds:
\begin{align}
\bm{\Sigma}\bm{\Sigma}^{\top} = \beta^{-1}\left(\bm{\Gamma} + \bm{\Gamma}^{\top}\right).
\label{EinsteinRelation}	
\end{align}
In our case, $\bm{\Sigma,\Gamma}$ are positive diagonal, and \eqref{EinsteinRelation} is equivalent to the proportionality constraint: $\sigma_{i}^{2} \propto \gamma_{i}$ for all $i=1,\hdots,n$.

In the power systems context, we cannot relate the damping coefficients $\gamma_i$ with the squared intensities of stochastic forcing $\sigma_{i}^{2}$ for the generators. Thus, \eqref{EinsteinRelation} will not hold in practice, meaning either we cannot guarantee existence-uniqueness for $\rho_{\infty}$, or even if $\rho_{\infty}$ exists, it will not be of the form \eqref{BoltzmannPDF}. On one hand, this implies that our construction of $\Phi$ may not be guided by free energy considerations. On the other hand, since we are only interested in computing the transient joint PDFs, i.e., \emph{non-equilibrium} statistical mechanics, the lack of a fluctuation-dissipation relation like \eqref{EinsteinRelation} should not be a fundamental impediment in setting up a recursion such as \eqref{ProxRecursionInfiniteDim}. We next show that a simple change of variable can indeed circumvent this issue.

\subsection{From Anisotropic to Isotropic Degenerate Diffusion}\label{subsec:FromAnisoToIso}
Consider the $2n\times 2n$ matrix
\begin{align}
\bm{\Psi}:= \bm{I}_{2} \otimes \left(\bm{M}\bm{\Sigma}^{-1}\right), \label{DefPsi}	
\end{align}
and define the invertible linear map
\begin{align}
\begin{pmatrix}
\bm{\theta}\\
\bm{\omega}	
\end{pmatrix} \mapsto \begin{pmatrix}
\bm{\xi}\\
\bm{\eta}	
\end{pmatrix} := \bm{\Psi} \begin{pmatrix}
\bm{\theta}\\
\bm{\omega}	
\end{pmatrix}. \label{ChangeOfVar}
\end{align}
Applying It\^{o}'s lemma \cite[Ch. 4.2]{oksendal2013stochastic} to the map \eqref{ChangeOfVar}, and using \eqref{ItoSDEvectorlevel}, we find that the transformed state vector $(\bm{\xi},\bm{\eta})^{\top}$ solves the It\^{o} SDE 
\begin{align}
\begin{pmatrix}
{\mathrm{d}}\bm{\xi}\\
{\mathrm{d}}\bm{\eta}	
\end{pmatrix}
\! = \!\begin{pmatrix}\!
\bm{\eta}\\
-\bm{\Upsilon}\nabla_{\bm{\xi}}U(\bm{\xi}) -\nabla_{\bm{\eta}}F\left(\bm{\eta}\right)  	
\!\end{pmatrix}\!{\mathrm{d}}t + \!\begin{pmatrix}\!
 \bm{0}_{n\times n}\\
 \bm{I}_{n}	
 \!\end{pmatrix}\!{\mathrm{d}}\bm{w},
\label{XiEtaVectorSDE}	
\end{align}
where the diagonal matrix $\bm{\Upsilon}:= \left(\prod_{i=1}^{n}\sigma_{i}^{2}m_{i}^{-2}\right)\bm{M\Sigma}^{-2}$ has positive entries along the main diagonal, and the potentials
{\small{\begin{subequations}
\begin{align}
U(\bm{\xi}) &:= \left(\displaystyle\sum_{i<j}k_{ij}\left(\!1 \!- \!\cos\!\left(\!\dfrac{\sigma_{i}}{m_{i}}\xi_{i} - \dfrac{\sigma_{j}}{m_{j}}\xi_{j} - \varphi_{ij}\!\!\right)\!\right) \right.\nonumber\\
&\left.\qquad\qquad-\displaystyle\sum_{i=1}^{n}\!\frac{\sigma_{i}}{m_{i}}P_{i}\xi_{i}\right)\left(\prod_{i=1}^{n}\frac{m_{i}^{2}}{\sigma_{i}^{2}}\right), \label{defU}
\\
F(\bm{\eta}) &:= \dfrac{1}{2}\langle\bm{\eta},\bm{\Sigma}^{-1}\bm{\Gamma}\bm{\eta}\rangle. \label{defF}	
\end{align}
\label{NewPotentials}	
\end{subequations}}}

Notice that \eqref{XiEtaVectorSDE} is a mixed conservative-dissipative SDE with \emph{isotropic} degenerate diffusion. In particular, the pushforward of the known initial joint PDF \eqref{FPKInitCond} via $\bm{\Psi}$, is given by
\begin{align}
\tilde{\rho}_{0}(\bm{\xi},\bm{\eta}) := \bm{\Psi}_{\sharp} \rho_{0} 
&= \dfrac{\rho_{0}\left(\bm{\Psi}^{-1} \begin{pmatrix}
\bm{\xi}\\
\bm{\eta}	
\end{pmatrix}\right)}{|\det\left(\bm{\Psi}\right)|} \nonumber\\
&= \dfrac{\rho_{0}\left(\bm{\Sigma}\bm{M}^{-1}\bm{\xi},\bm{\Sigma}\bm{M}^{-1}\bm{\eta}\right)}{\left(\displaystyle\prod_{i=1}^{n}m_{i}/\sigma_{i}\right)^{\!2}},
\label{PushforwardOfIC}	
\end{align}
where we used the standard properties of the Kronecker product. The transient joint state PDF $\tilde{\rho}(t,\bm{\xi},\bm{\eta})$ corresponding to \eqref{XiEtaVectorSDE} solves the PDE initial value problem
\begin{subequations}
\begin{align}
&\dfrac{\partial\tilde{\rho}}{\partial t} = -\langle\bm{\eta},\nabla_{\bm{\xi}}\tilde{\rho}\rangle + \nabla_{\bm{\eta}}\cdot\left(\tilde{\rho}\left(\bm{\Upsilon}\nabla_{\bm{
\xi}}U(\bm{\xi})+\nabla_{\bm{\eta}}F(\bm{\eta})\right)\right) \nonumber\\
&\qquad\qquad\qquad\qquad\qquad\qquad\qquad\qquad+ \dfrac{1}{2}\Delta_{\bm{\eta}}\tilde{\rho}, \label{NewPDE}\\
&\tilde{\rho}(t=t_{0},\bm{\xi},\bm{\eta}) = \underbrace{\tilde{\rho}_{0}(\bm{\xi},\bm{\eta})}_{\text{from}\;\eqref{PushforwardOfIC}}. \label{NewIC}	
\end{align}
\label{TransformedPDEivp}	
\end{subequations}
In other words, \eqref{TransformedPDEivp} is the macroscopic dynamics corresponding to the sample path dynamics \eqref{XiEtaVectorSDE}.

Since \eqref{NewPDE} is a kinetic Fokker-Planck PDE with isotropic degenerate diffusion, our strategy is to perform a proximal recursion of the form \eqref{ProxRecursionInfiniteDim} for \eqref{TransformedPDEivp} in $(\bm{\xi},\bm{\eta})$ coordinates, and then to pushforward the resulting joint PDFs via $\bm{\Psi}^{-1}$ to the original state space. This is what we detail next.

\begin{figure}[t]
\centering
\includegraphics[width=.7\linewidth]{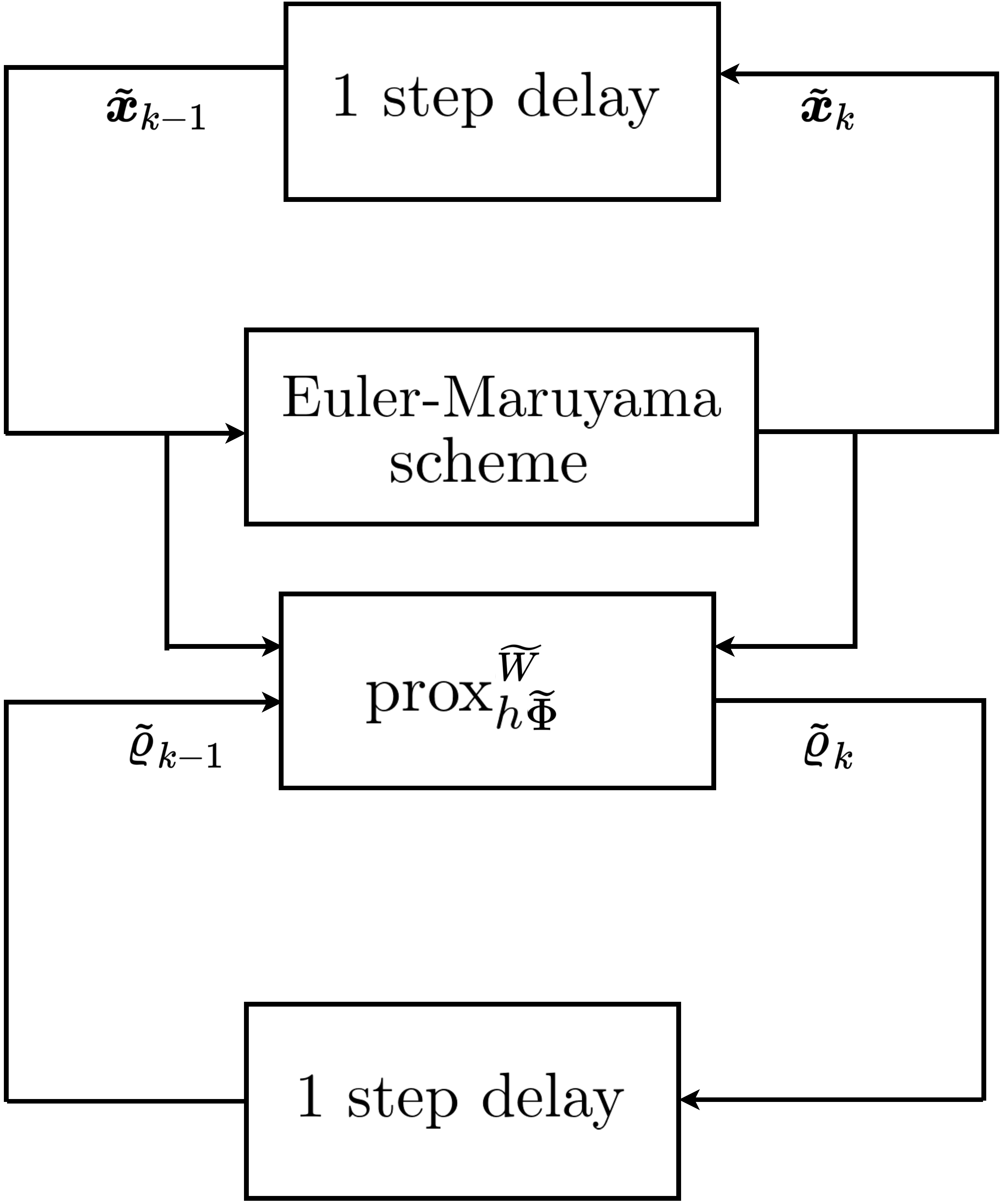}
\caption{\small{Schematic of the proposed proximal algorithm for propagating the joint state PDF as probability-weighted scattered point cloud $\{\tilde{\bm{x}}_{k}^{i},\tilde{\varrho}_{k}^{i}\}_{i=1}^{N}$. The states $\{\tilde{\bm{x}}_{k}^{i}\}_{i=1}^{N}$ are updated by the Euler-Maruyama scheme applied to (\ref{ItoSDEvectorlevel}); the corresponding probability weights $\{\tilde{\varrho}_{k}^{i}\}_{i=1}^{N}$ are updated via discrete version of the proximal recursion (\ref{KineticProx}) as detailed in Sec. \ref{sec:ProxAlgo}.}}
\label{fig:BlockDiagm}
\vspace*{-0.2in}
\end{figure}

\subsection{Proximal Update}\label{subsec:ProxUpdate}
Looking at \eqref{NewPotentials} and \eqref{NewPDE}, it is natural to consider an energy functional of the form
\begin{align}
\Phi(\tilde{\rho}) &:= \!\!\displaystyle\int_{\left(\prod_{i=1}^{n}[0,2\pi m_{i}/\sigma_{i})\right)\times\mathbb{R}^{n}}\!\!\!\left(\!U(\bm{\xi})+ F\left(\bm{\eta}\right) + \frac{1}{2}\log\tilde{\rho}\!\right)\nonumber\\
&\qquad\qquad\qquad\qquad\qquad\qquad\qquad\qquad\tilde{\rho}\:\differential\bm{\xi}\:\differential\bm{\eta},
\label{LyapButNotProx}	
\end{align}
which is the sum of a potential energy (expected value of $U$), a weighted kinetic energy (expected value of $F$), and an internal energy (scaled negative entropy, the entropy being $-\int\tilde{\rho}\log\tilde{\rho}$). 
However, unlike the gradient drift case mentioned in Sec. \ref{subsec:GenGradDescent}, it is not possible to express the right hand side of \eqref{NewPDE} as the Wasserstein gradient of the functional \eqref{LyapButNotProx}. This can be verified via direct computation by recalling that the Wasserstein gradient is defined as \cite[Ch. 8]{ambrosio2008gradient} 
$$\nabla^{\text{Wasserstein}}\Phi := -\nabla\cdot\left(\rho\nabla\dfrac{\delta\Phi}{\delta\rho}\right),$$
where $\nabla$ denotes the standard Euclidean gradient w.r.t. the vector $(\bm{\xi},\bm{\eta})^{\top}$, and $\frac{\delta}{\delta\rho}$ denotes the functional derivative. 
Thus, 
we cannot interpret the flow generated by \eqref{TransformedPDEivp} as the Wasserstein gradient flow of the functional \eqref{LyapButNotProx}. Consequently, in \eqref{InfiniteDimProxOpDef}, we cannot construct $(\Phi,{\rm{dist}})$ by pairing \eqref{LyapButNotProx} with the Wasserstein metric.

To set up a variational recursion of the form \eqref{InfiniteDimProxOpDef} for \eqref{TransformedPDEivp}, we set $(\Phi,{\rm{dist}}) \equiv \left(\widetilde{\Phi},\widetilde{W}_{h,\bm{\Upsilon}}\right)$ where
\begin{align}
\widetilde{\Phi}(\tilde{\rho}) := \!\!\int_{\left(\prod_{i=1}^{n}[0,2\pi m_{i}/\sigma_{i})\right) \times \mathbb{R}^{n}}\!\!\left(\!F(\boldsymbol{\eta})+\frac{1}{2} \log \tilde{\rho}\!\right) \tilde{\rho}\: \mathrm{d} \boldsymbol{\xi} \mathrm{d} \boldsymbol{\eta},
\label{Phitilde}	
\end{align}
and 
\begin{align}
&\widetilde{W}_{h,\bm{\Upsilon}}^{2}\left(\tilde{\varrho},\tilde{\varrho}_{k-1}\right) := \!\underset{\pi\in\Pi\left(\tilde{\varrho},\tilde{\varrho}_{k-1}\right)}{\inf}\!\displaystyle\int_{\left(\prod_{i=1}^{n}[0,2\pi m_{i}/\sigma_{i})\right)^{2} \times \mathbb{R}^{2n}}\nonumber\\
&\qquad\qquad\qquad\qquad s_{h,\bm{\Upsilon}}\left(\bm{\xi},\bm{\eta},\bar{\bm{\xi}},\bar{\bm{\eta}}\right)\differential\pi\left(\bm{\xi},\bm{\eta},\bar{\bm{\xi}},\bar{\bm{\eta}}\right),
\label{Wtilde}    
\end{align}
wherein $h>0$ denotes the step-size in proximal recursion, and $\Pi\left(\tilde{\varrho},\tilde{\varrho}_{k-1}\right)$ denotes the set of joint probability measures over the product space $\left(\prod_{i=1}^{n}[0,2\pi m_{i}/\sigma_{i})\right)^{2} \times \mathbb{R}^{2n}$ that have finite second moments  with the first marginal $\tilde{\varrho}$, and the second marginal $\tilde{\varrho}_{k-1}$. The ``ground cost" $s_{h,\bm{\Upsilon}}$ in \eqref{Wtilde} is given by
 \begin{align}
     &s_{h,\bm{\Upsilon}}\left(\bm{\xi},\bm{\eta},\bar{\bm{\xi}},\bar{\bm{\eta}}\right):=\nonumber\\
    &\bigg\langle\!\!\left(\bar{\bm{\eta}}-\bm{\eta}+h\bm{\Upsilon} \nabla_{\bm{\xi}} \bm{U}(\bm{\xi})\right),\bm{\Upsilon}^{-1}\!\left(\bar{\bm{\eta}}-\bm{\eta}+h\bm{\Upsilon} \nabla_{\bm{\xi}}\bm{U}(\bm{\xi})\right)\!\!\bigg\rangle\nonumber\\
    & + 12 \bigg\langle\!\left(\!\frac{\bar{\bm{\xi}}-\bm{\xi}}{h}-\frac{\bar{\bm{\eta}}-\bm{\eta}}{h}\!\right),\bm{\Upsilon}^{-1}\!\left(\!\frac{\bar{\bm{\xi}}-\bm{\xi}}{h}-\frac{\bar{\bm{\eta}}-\bm{\eta}}{h}\!\right)\!\bigg\rangle.
    \end{align}


That the sequence of functions $\{\tilde{\varrho}_{k}\}$ for $k\in\mathbb{N}$, generated by the proximal recursion
\begin{align}
\tilde{\varrho}_{k} &= {\text{prox}}^{\widetilde{W}}_{h\widetilde{\Phi}}(\tilde{\varrho}_{k-1}) \nonumber\\
&\equiv \underset{\tilde{\varrho}\in\mathcal{P}_{2}}{\arg\inf}\:\frac{1}{2}\widetilde{W}^{2}\left(\tilde{\varrho},\tilde{\varrho}_{k-1}\right) + h\:\widetilde{\Phi}(\tilde{\varrho}), \quad \tilde{\varrho}_{0} := \tilde{\rho}_{0},
\label{KineticProx}	
\end{align}
converges to the flow generated by \eqref{TransformedPDEivp}, i.e., 
$$\tilde{\varrho}_{k}(\bm{\xi},\bm{\eta}) \xrightarrow{h\downarrow 0} \tilde{\rho}(t=kh,\bm{\xi},\bm{\eta})\;\text{in}\;L^{1}\left(\mathbb{T}^{n}\times\mathbb{R}^{n}\right),$$
can be established following the arguments in \cite{duong2014conservative}. To numerically perform the recursion \eqref{KineticProx}, we employ the proximal algorithm proposed in \cite{caluya2019TAC} with finite number of samples, as explained next.


\section{Proximal Algorithm}\label{sec:ProxAlgo}
We solve \eqref{KineticProx} by recursively updating the probability-weighted scattered point clouds $\{\tilde{\bm{x}}_{k}^{i},\tilde{\varrho}_{k}^{i}\}_{i=1}^{N}$ where 
$$\tilde{\bm{x}}_{k}^{i}:=\left(\bm{\xi}_{k}^{i},\bm{\eta}_{k}^{i}\right)^{\!\top}, \quad i=1,\hdots,N, \quad k\in\mathbb{N}.$$
Thus, $\tilde{\varrho}_{k}^{i}$ is the joint PDF value obtained from \eqref{KineticProx} at $\tilde{\bm{x}}_{k}^{i}$, the $i$th (transformed) state sample at the $k$th time step. The high level schematic of the algorithm is shown in Fig. \ref{fig:BlockDiagm}.

In the numerical simulations reported in Sec. \ref{sec:NumericalSimulations}, the states $\{\tilde{\bm{x}}_{k}^{i}\}_{i=1}^{N}$ are updated by the Euler-Maruyama scheme applied to (\ref{ItoSDEvectorlevel}). If one wishes so, the Euler-Maruyama scheme in Fig. \ref{fig:BlockDiagm} may be replaced by other SDE integrators, see e.g., \cite[Sec. III.B.2, Remark 1]{caluya2019TAC}. 

To numerically perform the proximal updates $\{\tilde{\varrho}_{k-1}^{i}\}_{i=1}^{N}\mapsto\{\tilde{\varrho}_{k}^{i}\}_{i=1}^{N}$ for $k\in\mathbb{N}$, we implement an instance of the Algorithm 1 in \cite{caluya2019TAC}. The algorithm involves a dualization along
with an entropic regularization of the variational update \eqref{KineticProx}, and then solving the same using a fixed point recursion that is provably contractive; we refer the interested readers to \cite[Sec. V.B]{caluya2019TAC} for details. This enables a nonparametric computation of $\tilde{\varrho}_{k}^{i} \equiv \tilde{\varrho}_{k}\left(\bm{\xi}_{k}^{i},\bm{\eta}_{k}^{i}\right)$ for $i=1,\hdots,N$.

\begin{algorithm}[t]
\caption{Proposed proximal algorithm for $\tilde{\bm{\varrho}}_{k-1}\mapsto\tilde{\bm{\varrho}}_{k}$}
\label{algo:KineticProx}
\begin{algorithmic}[1]
\Procedure{Prox}{$\tilde{\bm{\varrho}}_{k-1}, \tilde{\bm{x}}_{k-1}, \tilde{\bm{x}}_{k}, \bm{\Upsilon}, h, \varepsilon, N, \delta, \ell_{\max}$}
\For{$i=1$ to $N$}
\State $\bm{\zeta}_{k-1}(i) \gets \exp\left(-F\left(\bm{\eta}_{k-1}^{i}\right)-1\right)$
\For{$j=1$ to $N$}
\State $\bm{C}_{k}(i,j) \gets s_{h,\bm{\Upsilon}}\left(\bm{\xi}_{k-1}^{i},\bm{\eta}_{k-1}^{i},\bm{\xi}_{k}^{j},\bm{\eta}_{k}^{j}\right)$
\EndFor
\EndFor
\State $\bm{\Gamma}_{k} \gets \exp\left(-\bm{C}_{k}/2\varepsilon\right)$ \Comment{elementwise exponential}
\State $\bm{z}_{0} \gets {\rm{rand}}_{N\times 1}$ \Comment{random vector of size $N\times 1$}
\State $\bm{z} \gets \left[\bm{z}_{0}, \bm{0}_{N\times(L-1)}\right]$ \Comment{initialize}
\State $\bm{y} \gets \left[\tilde{\bm{\varrho}}_{k-1}\oslash\left(\bm{\Gamma}_{k}\bm{z}_{0}\right), \bm{0}_{N\times(L-1)}\right]$ \Comment{initialize}
\State $\ell = 1$
\While{$\ell\leq \ell_{\max}$}
\State $\bm{z}(:,\ell+1) \gets \left(\bm{\zeta}_{k-1}\oslash\left(\bm{\Gamma}_{k}^{\top}\bm{y}(:,\ell)\right)\right)^{\frac{1}{1+2\varepsilon/h}}$
\State $\bm{y}(:,\ell+1) \gets \tilde{\bm{\varrho}}_{k-1} \oslash \left(\bm{\Gamma}_{k}\bm{z}(:,\ell+1)\right)$
\If{$\|\bm{y}(:,\ell+1)-\bm{y}(:,\ell)\|_{2}<\delta$ \& $\|\bm{z}(:,\ell+1)-\bm{z}(:,\ell)\|_{2}<\delta$} \Comment{error within tolerance}
\State break
\Else
\State $\ell \gets \ell + 1$
\EndIf
\EndWhile\\
\Return $\tilde{\bm{\varrho}}_{k} \gets \bm{z}(:,\ell) \odot \left(\bm{\Gamma}_{k}^{\top}\bm{y}(:,\ell)\right)$ \Comment{proximal update}
\EndProcedure	
\end{algorithmic}
\end{algorithm}

Finally, we transform the proximal updates back to the $\bm{x}\equiv\left(\bm{\theta},\bm{\omega}\right)^{\top}$ state space via the pushforward $\bm{\Psi}^{-1}$ as
\begin{multline}
\varrho_{k}\left(\bm{\theta}_{k}^{i},\bm{\omega}_{k}^{i}\right) = \bm{\Psi}^{-1}_{\sharp} \tilde{\varrho}_{k}\left(\bm{\xi}_{k}^{i},\bm{\eta}_{k}^{i}\right) \\
= \left(\prod_{j=1}^{n}m_{j}/\sigma_{j}\right)^{\!\!2}\tilde{\varrho}_{k}\left(\bm{M}\bm{\Sigma}^{-1}\bm{\theta}_{k}^{i},\bm{M}\bm{\Sigma}^{-1}\bm{\omega}_{k}^{i}\right),
\label{BringBackOriginalStateSpace}	
\end{multline}
for all $i=1,\hdots,N$. The $\varrho_{k}\left(\bm{\theta},\bm{\omega}\right)$ from \eqref{BringBackOriginalStateSpace} approximates $\rho(t,\bm{\theta},\bm{\omega})$ (the transient solution of \eqref{KineticFPK}) in the sense \eqref{ConvergenceGuarantee}. 

For completeness, the algorithm \textproc{Prox} for updating $\{\tilde{\varrho}_{k-1}^{i}\}_{i=1}^{N}\mapsto\{\tilde{\varrho}_{k}^{i}\}_{i=1}^{N}$ is outlined in Algorithm \ref{algo:KineticProx}. As shown in Fig. \ref{fig:BlockDiagm}, this algorithm, at a conceptual level, takes the pre and post-update state samples 
\[\{\tilde{\bm{x}}_{k-1}^{i}\}_{i=1}^{N}\equiv\{\left(\bm{\xi}_{k-1}^{i},\bm{\eta}_{k-1}^{i}\right)\}_{i=1}^{N}, \; \{\tilde{\bm{x}}_{k}^{i}\}_{i=1}^{N}\equiv\{\left(\bm{\xi}_{k}^{i},\bm{\eta}_{k}^{i}\right)\}_{i=1}^{N},\]
and $\{\tilde{\varrho}_{k-1}^{i}\}_{i=1}^{N}$ as inputs, and outputs the proximal updates $\{\tilde{\varrho}_{k}^{i}\}_{i=1}^{N}$. For each $k\in\mathbb{N}$, the updated probability-weighted point clouds $\{\tilde{\bm{x}}_{k}^{i},\tilde{\varrho}_{k}^{i}\}_{i=1}^{N}$ are then brought back to the original state space via $\bm{\Psi}^{-1}$ as $\{\bm{x}_{k}^{i},\varrho_{k}^{i}\}_{i=1}^{N}$, as explained earlier.

Algorithm \ref{algo:KineticProx} also needs input parameters $h,\varepsilon,N,\delta,L$. Specifically, $h$ is the time-step size, $\varepsilon$ is an entropic regularization weight internal to the \textproc{Prox} algorithm, and $N$ is the number of samples. The parameters $\delta$ and $\ell_{\max}$ codify the numerical tolerance and maximum number of iterations, respectively, for the while loop in Algorithm \ref{algo:KineticProx}. Its convergence guarantees can be found in \cite[Sec. III.C]{caluya2019TAC}.

\begin{remark}\label{remark:timevaryingEffectivePowerInput}
The proposed computational framework is also applicable when $P_{i}^{\rm{mech}}, P_{i}^{\rm{load}}$, and thus $P_{i}$ in \eqref{defPi}, are bounded time-varying functions. This can be the case, for instance, with rapid fluctuations from renewables or from loads in the short term.
\end{remark}


\section{Numerical Simulations}\label{sec:NumericalSimulations}
To illustrate the proposed computational framework, we next provide two numerical simulation case studies. In Sec. \ref{subsec:NumSimIEEEBus}, we consider the prediction of transient stochastic states for the IEEE 14 bus system for the nominal case as well as for the case when a line failure occurs. In Sec. \ref{subsec:NumSimSynthetic}, we propagate the joint PDFs over the 100 dimensional state space of a synthetic power network with randomly generated parameters--our intent there is to highlight the scalability of the proposed method. All simulations were performed in MATLAB R2019b on an iMac with 3.4 GHz Quad-Core Intel Core i5 processor and 8 GB memory.

\begin{figure*}[t]
\centering
\includegraphics[width=0.9\linewidth]{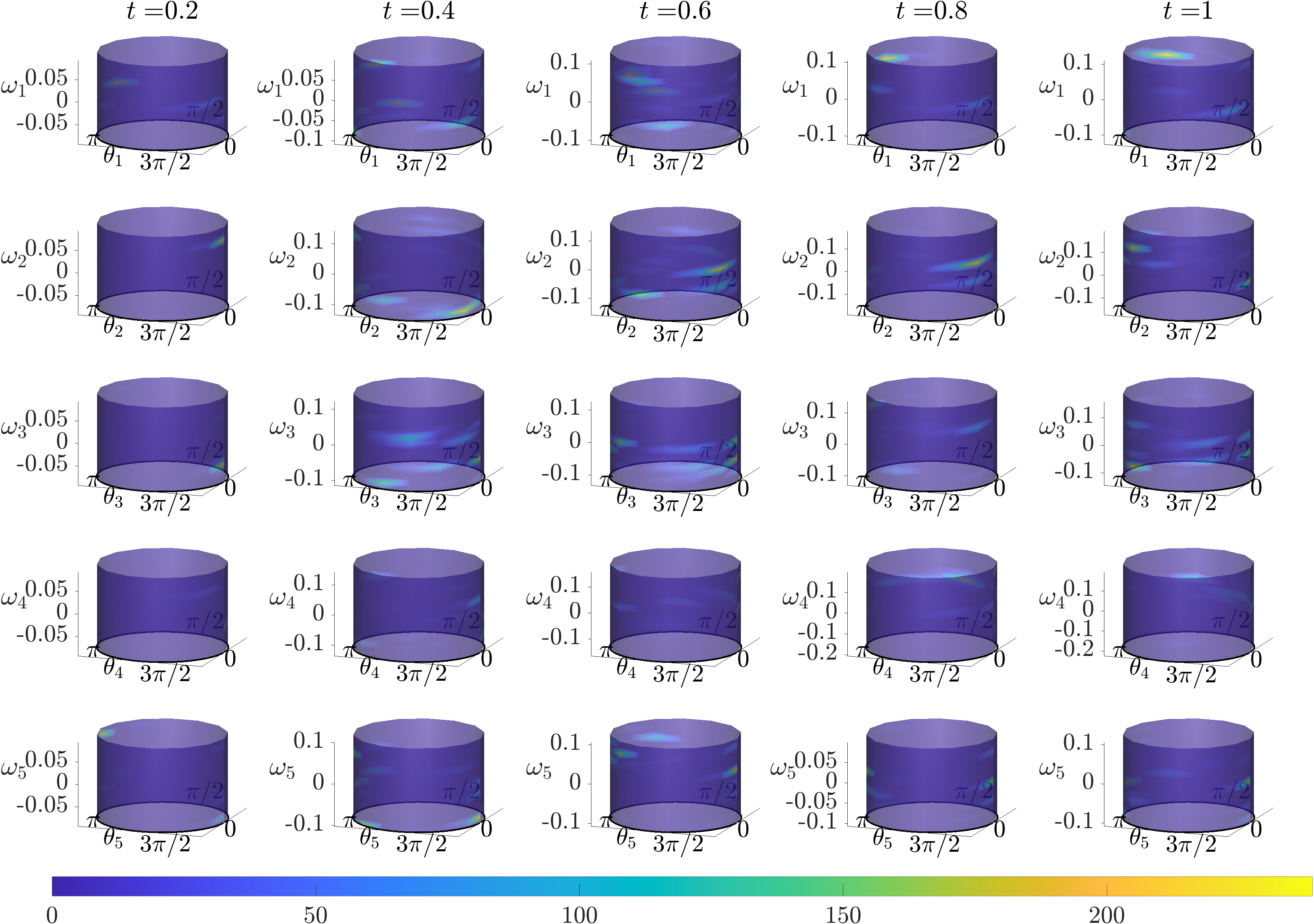}
\caption{\small{Time evolution of the $(\theta_{i},\omega_{i})$ bivariate marginal PDFs of the generator $i=1,\hdots,5$, for the Kron-reduced IEEE 14 bus simulation setup \textbf{Case I} described in Sec. \ref{subsec:NumSimIEEEBus}. The five rows above correspond to the five generator nodes; these are buses 1, 2, 3, 6 and 8 in Fig. \ref{fig:Motivating}. The columns correspond to the time snapshots. The colorbar at the bottom shows the marginal PDF values. The $\theta_i$ are in rad, the $\omega_i$ are in rad/s.}}
\vspace*{-0.1in}
\label{fig:IEEE14marginalsCase0}
\end{figure*}

\begin{figure*}[t]
\centering
\includegraphics[width=0.9\linewidth]{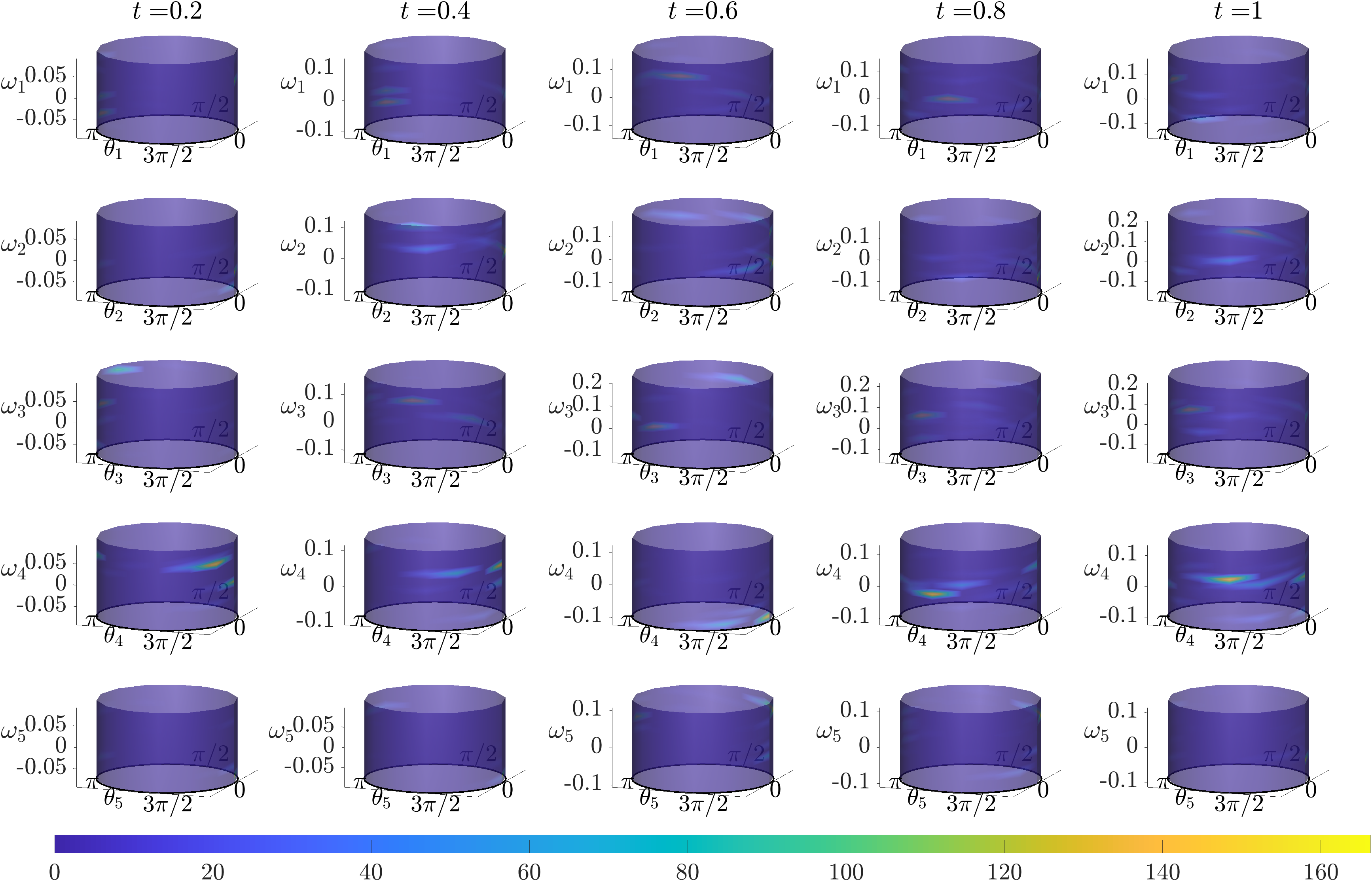}
\caption{\small{Time evolution of the $(\theta_{i},\omega_{i})$ bivariate marginal PDFs of the generator $i=1,\hdots,5$, for the Kron-reduced IEEE 14 bus simulation setup \textbf{Case II} described in Sec. \ref{subsec:NumSimIEEEBus}. The five rows above correspond to the five generator nodes; these are buses 1, 2, 3, 6 and 8 in Fig. \ref{fig:Motivating}. The columns correspond to the time snapshots. The colorbar at the bottom shows the marginal PDF values. The $\theta_i$ are in rad, the $\omega_i$ are in rad/s.}}
\vspace*{-0.1in}
\label{fig:IEEE14marginalsCase1}
\end{figure*}

\begin{figure*}[ht!]
   \subfigure[]{\includegraphics[width=0.322\textwidth]{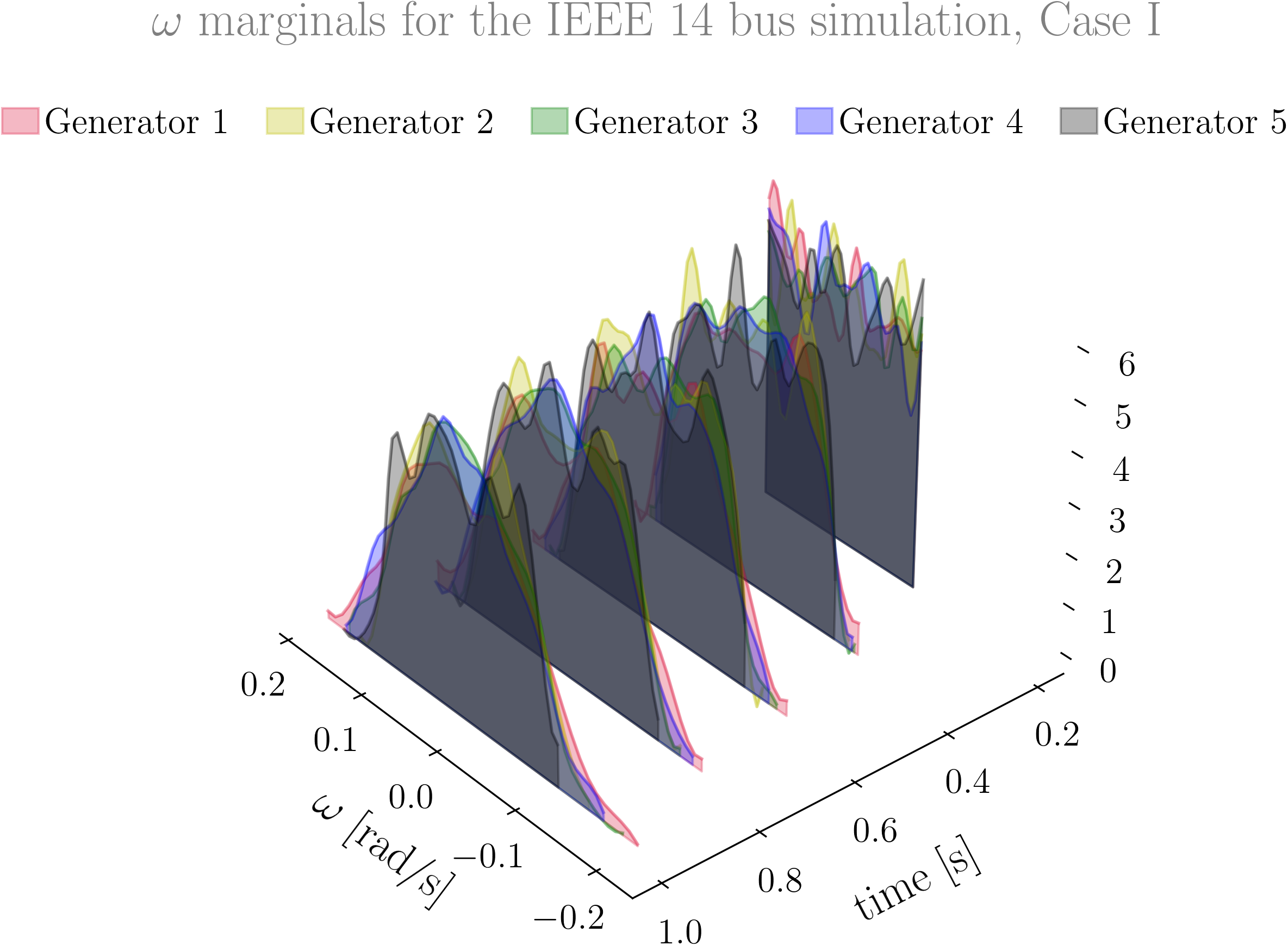}}
   \hspace{\fill}
   \subfigure[]{\includegraphics[width=0.322\textwidth]{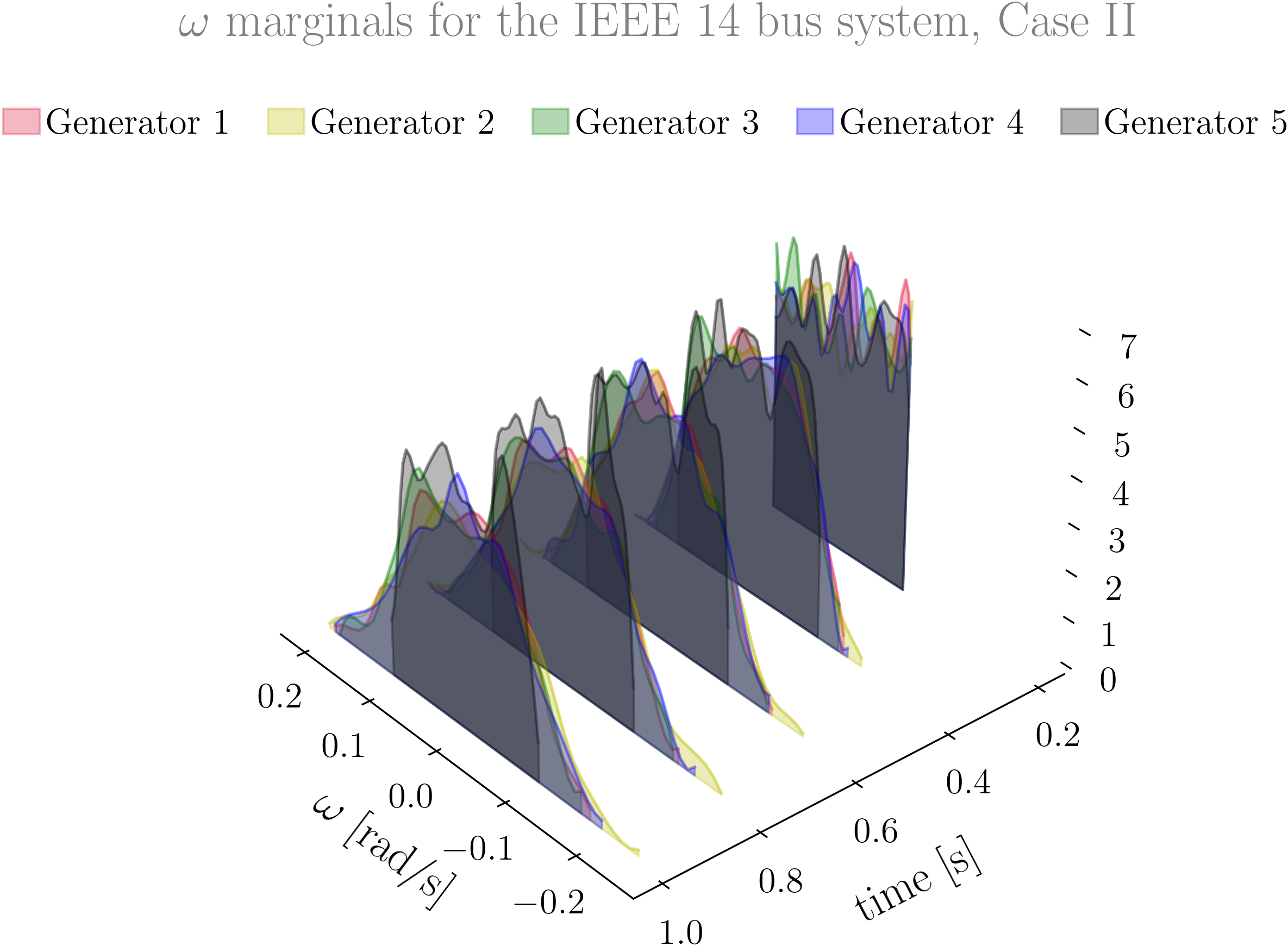}}
   \hspace{\fill}
   \subfigure[]{\includegraphics[width=0.309\textwidth]{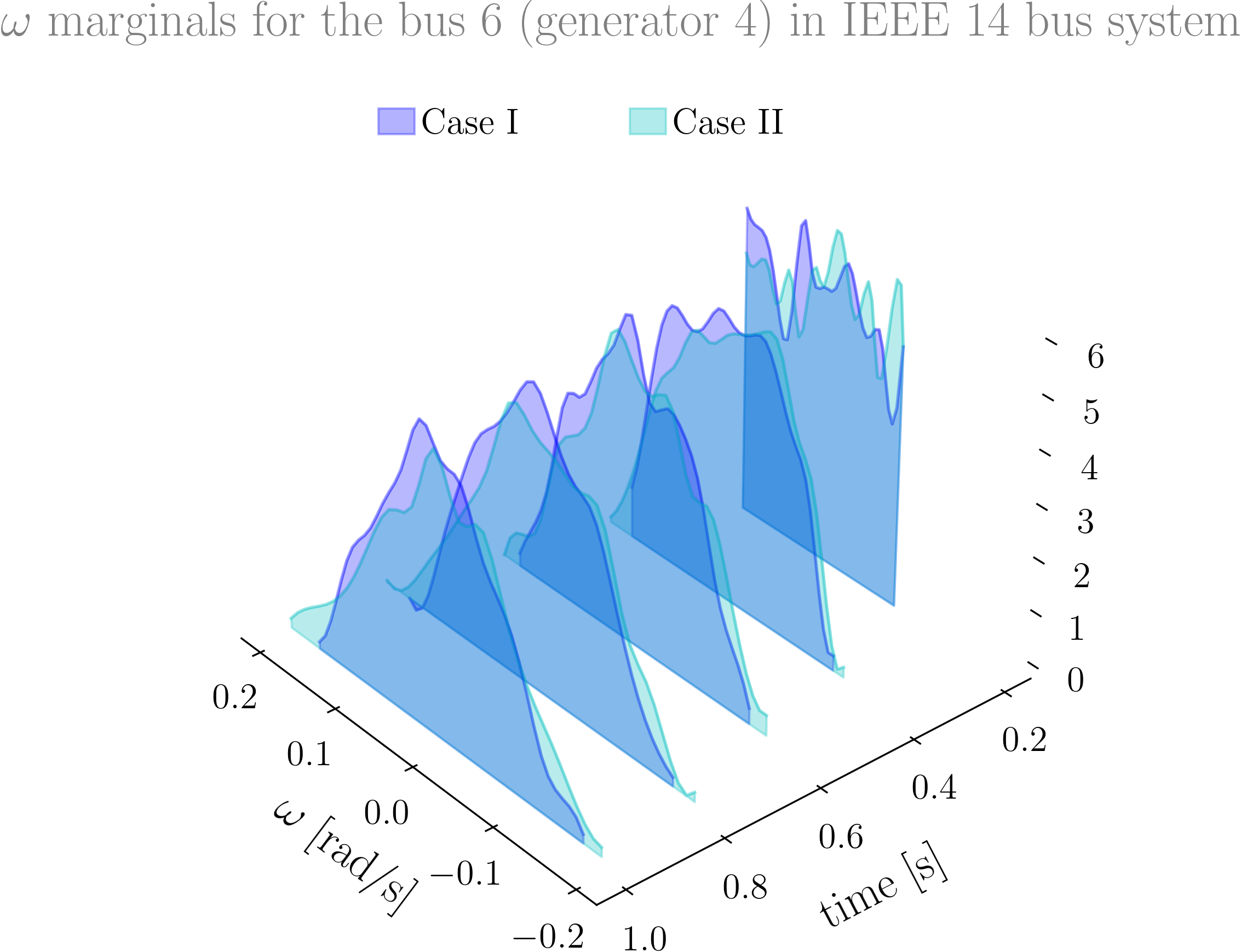}}
\caption{\small{For $i=1,\hdots,5$, the $\omega_{i}$ marginals (a) for \textbf{Case I}, and (b) for \textbf{Case II}, resulting from the IEEE 14 bus simulation described in Sec. \ref{subsec:NumSimIEEEBus}. (c) Comparison of the $\omega_{4}$ marginals for bus 6 in the nominal (\textbf{Case I}) and post-fault (\textbf{Case II}) scenarios.}}
\vspace*{-0.1in}
\label{fig:IEEE14omegamarginalsCase12}
\end{figure*}

\begin{figure*}[ht!]
   \subfigure[Case I]{\includegraphics[width=0.49\textwidth]{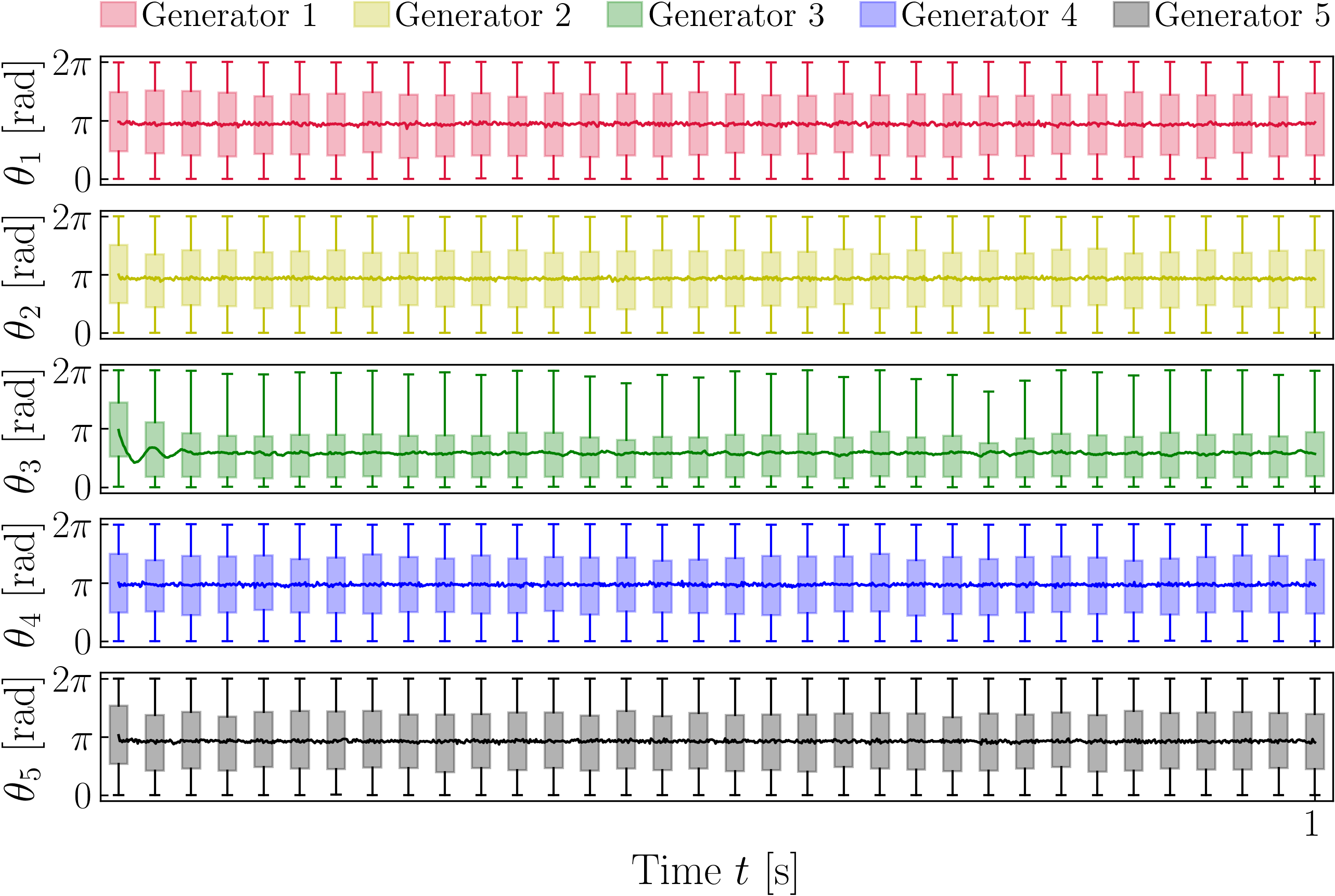}}
   \hspace{\fill}
   \subfigure[Case II]{\includegraphics[width=0.49\textwidth]{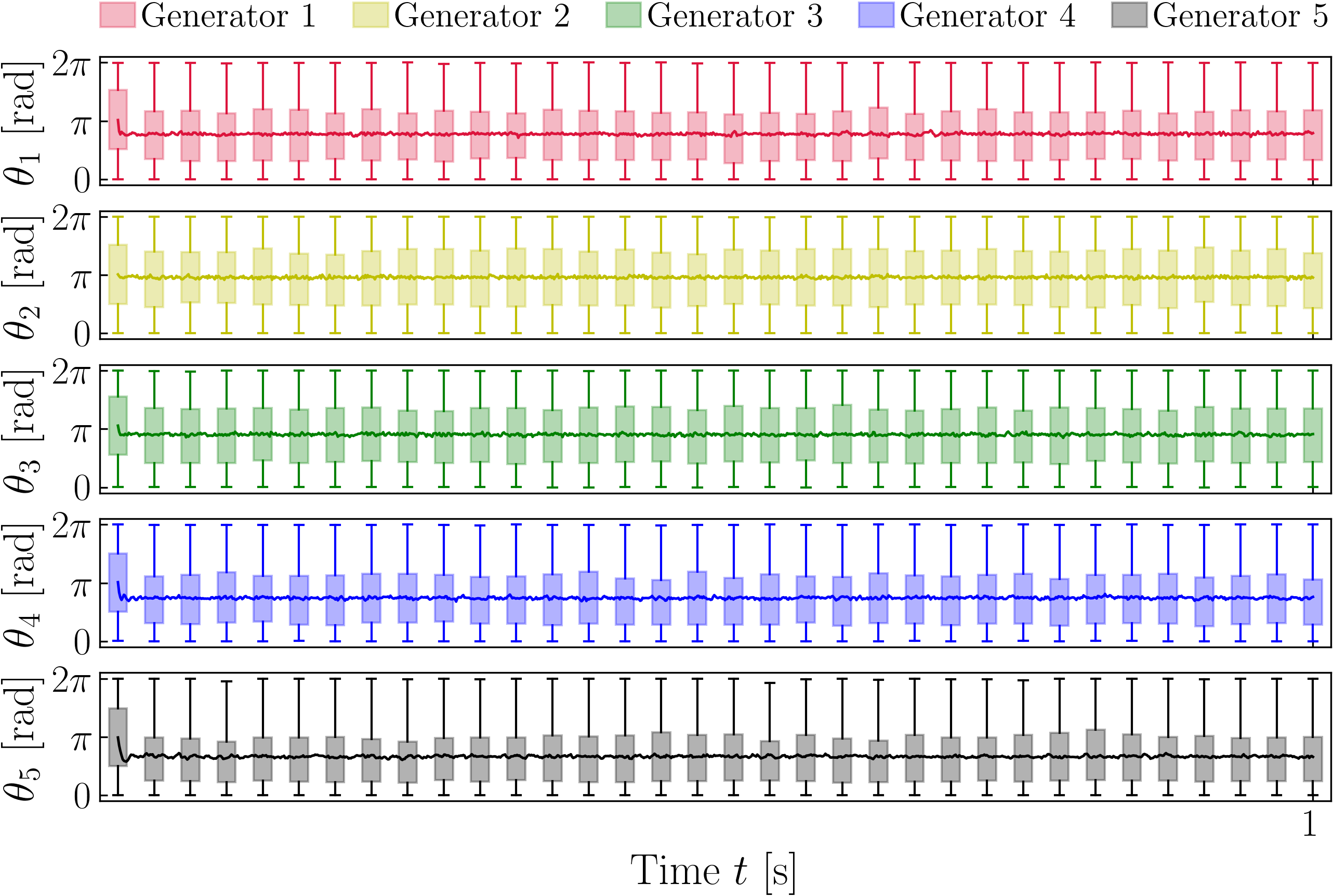}}
\caption{\small{Boxplots for the rotor angles $\theta_{i}(t)$, $i=1,\hdots,5$, for the IEEE 14 bus system simulation in Sec. \ref{subsec:NumSimIEEEBus} in the (a) nominal \textbf{Case I}, and in the (b) post-fault \textbf{Case II} scenarios, respectively. As usual, the boxes show the respective interquartile ranges, and the whiskers show the respective minimum and maximum values. In each box plot, we have superimposed  the respective \emph{mean} rotor angle sample path.}}
\vspace*{-0.1in}
\label{fig:IEEE14theta}
\end{figure*}

\begin{figure}[htpb]
\centering
\includegraphics[width=0.95\linewidth]{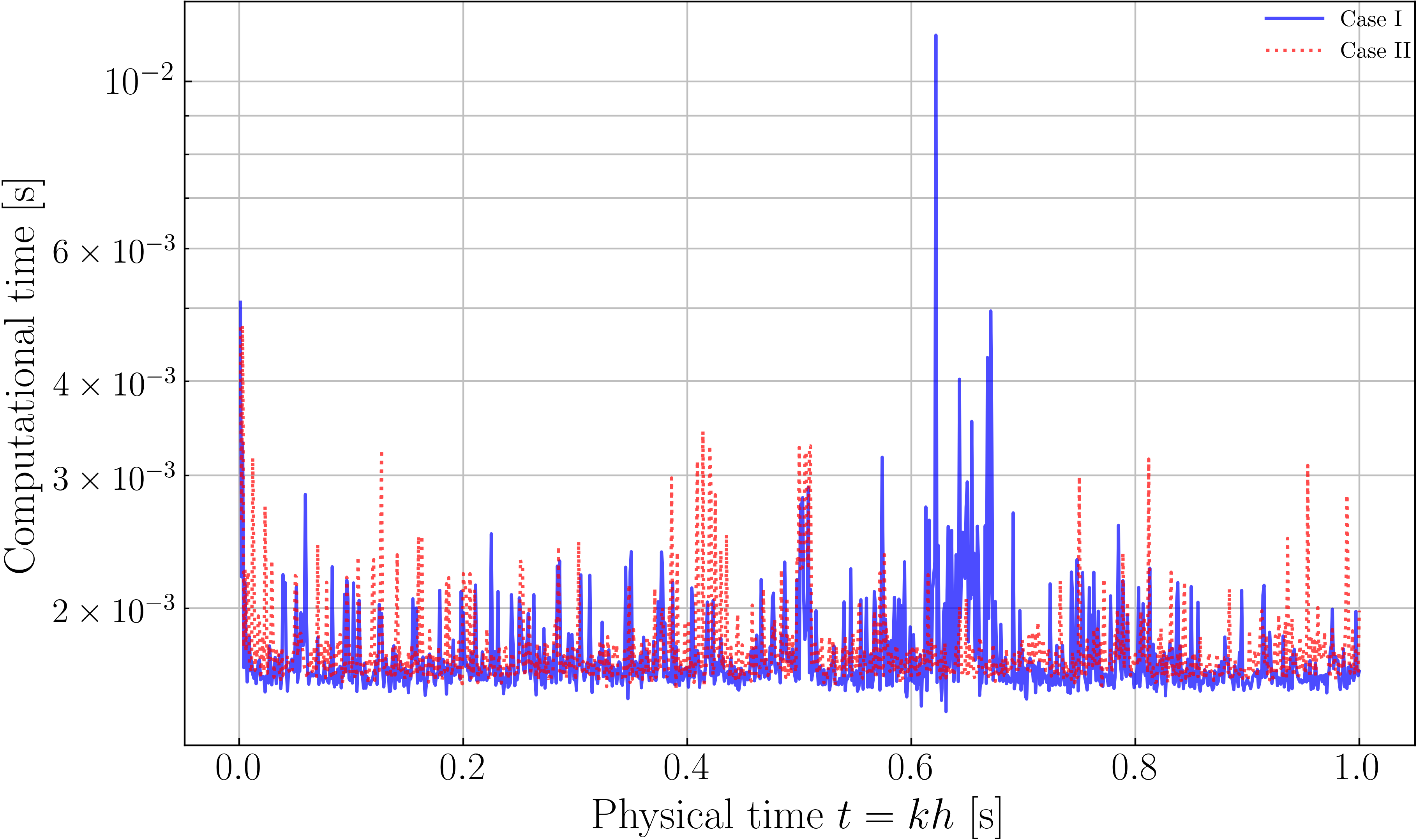}
\caption{\small{The computational times for propagating the transient joint state PDFs over $\mathbb{T}^{5}\times\mathbb{R}^{5}$ for the simulation set up in Sec. \ref{subsec:NumSimIEEEBus}.}}
\vspace*{-0.1in}
\label{fig:CompTimeIEEE14bus}
\end{figure}

\subsection{IEEE 14 Bus System}\label{subsec:NumSimIEEEBus}
We consider the Kron-reduced dynamics \eqref{ItoSDEcomponentlevel} for the IEEE 14 bus system shown in Fig. \ref{fig:Motivating}. In this case, we have $n=5$ nodes which correspond to the buses 1, 2, 3, 6 and 8 in Fig. \ref{fig:Motivating}. We obtained the parameters of the 14-bus system from MATPOWER \cite{zimmerman2010matpower}. The calculation of Kron-reduced admittance matrix and current vector in \eqref{DefParam} followed Sec.  \ref{subsec:samplepathdynamics}; see also \cite{dorfler2012kron}. The parameters $m_{i}, \gamma_{i}$ were obtained from the open-source Python-based
library ANDES \cite{cui2021andes}.

Our numerical simulations for the IEEE 14 bus system consider the following two cases:

\noindent\textbf{Case I.} \textbf{Nominal case} where the $P_{i}^{\text{mech}}$ and $P_{i}^{\text{load}}$ are from the steady state power flow solutions without the process noise, and are assumed to remain constant over the simulation time horizon. The transient stochastic simulations (i.e., with process noise) are then performed using the proposed framework.

\noindent\textbf{Case II.} \textbf{Post-contingency case} where the line 13 (see Fig. \ref{fig:Motivating}) fails at time $t=0$, and starting with the parameters and initial conditions of \textbf{Case I}, we simulate the post-contingency stochastic state evolution with process noise. 


For both the cases mentioned above, we used the randomly generated noise coefficients $\sigma_{i}\in [1,5]$ as
\[(\sigma_{1},\hdots,\sigma_5) = \left(2.4628, 4.9266, 4.8724, 1.4215, 3.8681\right).\]

To account the pertinent geometry of $\mathbb{T}^{n}\times\mathbb{R}^{n}$, we take the initial joint PDF as
\begin{align}
\rho_{0}&\equiv\rho(t=0,\bm{\theta}(0),\bm{\omega}(0)) = \vartheta_{0}\left(\bm{\theta}(0)\right)\Omega_{0}\left(\bm{\omega}(0)\right)\nonumber\\
&=\underbrace{\left(\prod_{i=1}^{n=5}\dfrac{\exp\left(\kappa_{i}\cos\left(2\theta_{i}(0)-\mu_{\theta_{i}(0)}\right)\right)}{2\pi I_{0}(\kappa_{i})}\right)}_{=:\vartheta_{0}\left(\bm{\theta}(0)\right)}\nonumber\\
	&\qquad\times \underbrace{{\rm{Unif}}\left(\left([-0.1,0.1]\;\text{rad/s}\right)^{5}\right)}_{=:\Omega_{0}\left(\bm{\omega}(0)\right)},
\label{rho0IEEE14}		
\end{align}
whose $\bm{\theta}$ marginal $ \vartheta_{0}\left(\bm{\theta}(0)\right)$ is a product von Mises PDF\footnote{the factor $2$ inside the cosine ensures the range $[0,2\pi)$; it can be dispensed if we instead use the range $[-\pi,\pi)$ for the angular variables, as common in the directional statistics literature \cite{mardia2009directional}.} \cite{mardia2008multivariate,mardia2014some} supported over $\mathbb{T}^{5}$ with mean angles (in rad)
\begin{align}
\left(\mu_{\theta_{1}}(0),\hdots,\mu_{\theta_{5}}(0)\right)=\left(0, 6.1963, 6.0612, 6.0350, 6.0500\right),
\label{MeanAngleInitialIEEE14}	
\end{align}
and concentration parameter $\kappa_{i}\geq 0$ set as 
\[(\kappa_{1},\hdots,\kappa_5) = \left(5,6,7,4,5\right).\]
In \eqref{rho0IEEE14}, $I_{0}(\cdot)$ denotes the modified Bessel function of the first kind with order zero, given by $I_{0}(\kappa) = \sum_{r=0}^{\infty} \frac{\left(\kappa^{2}/4\right)^{r}}{(r!)^{2}}$, and can be evaluated via MATLAB command \texttt{besseli(0,$\cdot$)}. 

The mean angles \eqref{MeanAngleInitialIEEE14} were obtained from steady state AC power flow solutions. Larger $\kappa_{i}$ entails a higher concentration around the mean angle $\mu_{\theta_{i}}(0)$ (setting $\kappa_i=0$ implies that $\theta_{i}(0)$ follows uniform distribution over $[0,2\pi)$).

In \eqref{rho0IEEE14}, the $\bm{\omega}$ marginal $ \Omega_{0}\left(\bm{\omega}(0)\right)$ is a uniform PDF over the hypercube $\left([-0.1,0.1]\;\text{rad/s}\right)^{5}$. Thus, the PDF $\rho_0$ in \eqref{rho0IEEE14} is supported over $\mathbb{T}^{5}\times [-0.1,0.1]^{5}$.

With step size $h=10^{-3}$, we discretized time $t\in[0,1]$. Using the algorithm detailed in Sec. \ref{sec:ProxAlgo}, for both cases, we propagated $N=1000$ random state samples from \eqref{rho0IEEE14}, and corresponding joint PDF values evaluated at \eqref{rho0IEEE14}, to update the scattered weighted point clouds $\{\bm{x}_{k}^{i},\varrho_{k}^{i}\}_{i=1}^{N}$ over $\mathbb{T}^{5}\times\mathbb{R}^{5}$ at times $t_{k}:=kh$, $k\in\mathbb{N}$. In Algorithm \ref{algo:KineticProx}, we used the algorithmic parameters $\varepsilon=0.05$, $\delta=10^{-3}$, $\ell_{\max}=300$.

Fig. \ref{fig:IEEE14marginalsCase0} and Fig. \ref{fig:IEEE14marginalsCase1} show the time snapshots of the $(\theta_{i},\omega_{i})$ bivariate marginal PDFs for \textbf{Case I} and \textbf{Case II} respectively, associated with the generator $i=1,\hdots, n=5$. These marginals are computed using the corresponding weighted point clouds $\{\bm{x}_{k}^{i},\varrho_{k}^{i}\}_{i=1}^{N}$.


Fig. \ref{fig:IEEE14omegamarginalsCase12} shows the corresponding univariate $\omega_{i}$-marginal PDFs. In particular, Fig. \ref{fig:IEEE14omegamarginalsCase12}(c) helps compare the $\omega_{4}$ marginal evolution for bus 6 in the nominal (\textbf{Case I}) and post-fault (\textbf{Case II}) scenarios. Since \textbf{Case II} considers the failure of line 13, from Fig. \ref{fig:Motivating} we intuitively expect that the transient effects of this fault will be prominent for bus 6 (generator 4). Indeed, Fig. \ref{fig:IEEE14omegamarginalsCase12}(c) shows that although the $\omega_{4}$ marginals in both \textbf{Case I} and \textbf{II}, at all times have the highest probability around $0$ rad/s, the post-fault marginals show larger spread or dispersion of the probability mass, as expected. That the generator 4 states in \textbf{Case II} have larger uncertainties than in \textbf{Case I}, is also visible in the bivariate marginal plots if we compare the fourth row in Fig. \ref{fig:IEEE14marginalsCase0} with the same in Fig. \ref{fig:IEEE14marginalsCase1}. The statistical uncertainties in the rotor angles are shown in Fig. \ref{fig:IEEE14theta} for both \textbf{Case I} and \textbf{Case II}.

Fig. \ref{fig:CompTimeIEEE14bus} shows the computational times needed to update $\{\bm{x}_{k}^{i},\varrho_{k}^{i}\}_{i=1}^{N}$ via the algorithm in Sec. \ref{sec:ProxAlgo}, for both \textbf{Cases I} and \textbf{II}. While the depiction of the computational times in Fig. \ref{fig:CompTimeIEEE14bus} appears unconventional, we next explain why that is meaningful in our framework.

At each fixed physical time $t_{k}=kh$, $k\in\mathbb{N}$, where $h$ is the fixed step-size, we perform the proximal update $\varrho_{k-1}\mapsto\varrho_{k}$ via Algorithm 1. This proximal update has no analytical solution and instead requires execution of a contractive fixed point recursion. The physical time $t_{k}$ is ``frozen" (zero-order-hold) during the fixed point recursions underlying this one-step proximal update. For the proposed method to be practical, we therefore, need to demonstrate that the computational timescale is at least \emph{statistically} comparable (in the order of magnitude sense) to the physical dynamics timescale. In other words, the computational times needed to complete the zero-order-hold fixed point recursions are, with high probability, the same orders-of-magnitude as $h$, the step-size of physical time. The precise computational time is probabilistic, not just due to random sampling, but also due to randomness in latency, power and memory footprint of the hardware platform.


In Fig \ref{fig:CompTimeIEEE14bus}, $h=10^{-3}$ and the computational time needed to perform a single proximal update, with high probability, is $O(10^{-3})$ in the sense most time-steps incur $O(10^{-3})$ time. For the Case I in Fig \ref{fig:CompTimeIEEE14bus}, few instances when computational times are $O(10^{-2})$ do occur but are statistically rare.

\begin{figure}[htpb]
\centering
\includegraphics[width=0.95\linewidth]{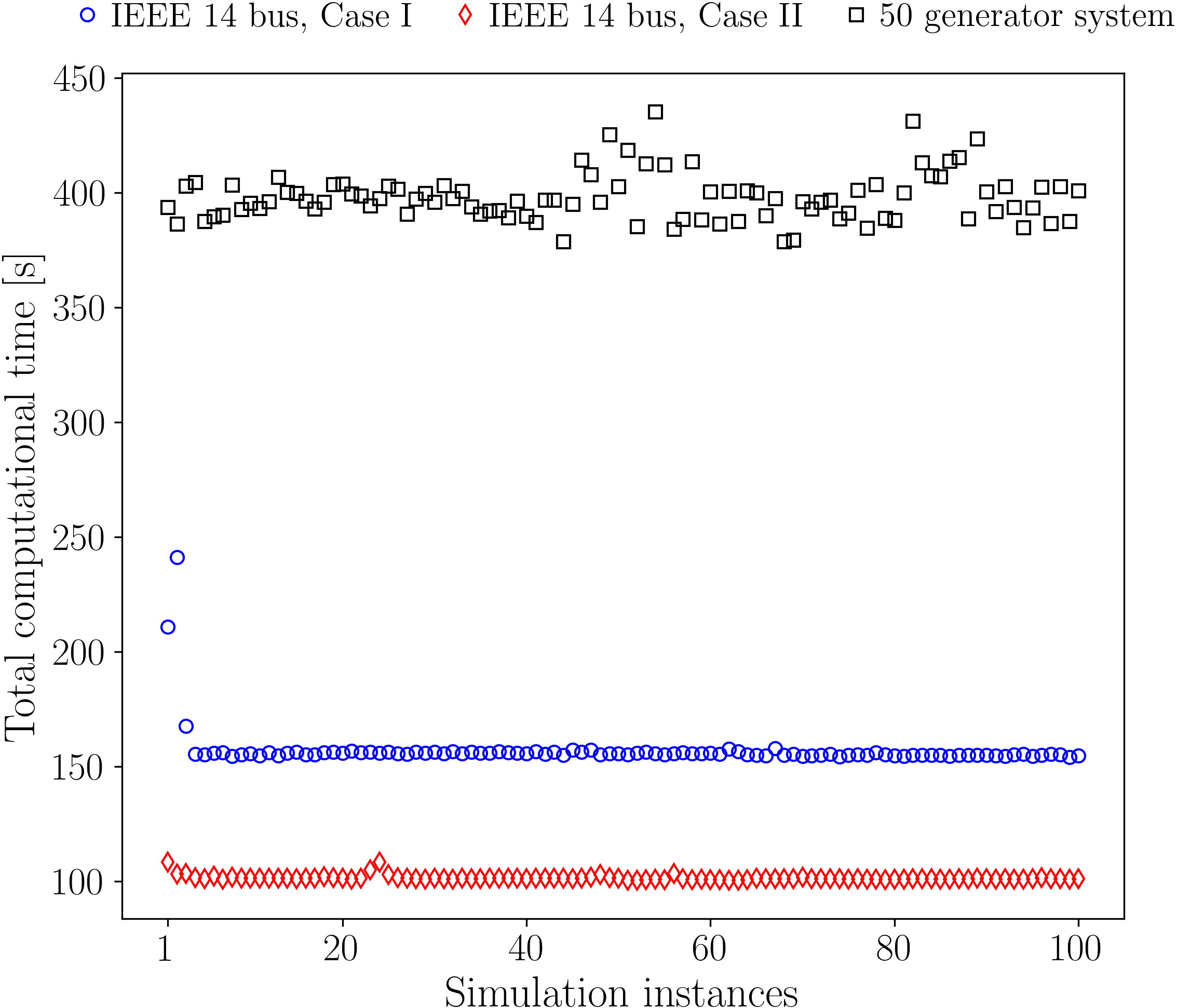}
\caption{\small{The total computational times for propagating the transient joint state PDFs over $[0, 1\;\text{min}]$ for the simulation set up in Sec. \ref{subsec:NumSimIEEEBus} and Sec. \ref{subsec:NumSimSynthetic}. In each case, we show the total computational time for 100 different simulation instances.}}
\vspace*{-0.1in}
\label{fig:TotalCompTime}
\end{figure}

Fig. \ref{fig:TotalCompTime} (circle markers for the \textbf{Case I}, diamond markers for the \textbf{Case II}) shows the total computational times needed for propagating the transient joint state PDFs over a longer horizon $[0, 1\;\text{min}]$ for 100 different simulation instances.

\begin{table}[t]
\centering
\begin{tabular}{| l | l |} 
\hline
Parameter description & Values\\
\hline\hline
&\\
nominal frequency & $f_{0} = 60$ Hz\\
&\\
inertia & $m_{i}\in{\rm{Unif}}\left([2,12]\right)/2\pi f_{0}$\\
& \\
damping coefficient & $\gamma_{i}\in{\rm{Unif}}\left([20,30]\right)/2\pi f_{0}$\\
& \\
diffusion coefficient & $\sigma_{i}\in{\rm{Unif}}\left([1,5]\right)$\\
& \\
tangent of phase shift & $\tan\varphi_{ij}\begin{cases}
=0 & \text{for}\; i=j\\
\in{\rm{Unif}}\left([0,0.25]\right)& \text{for}\; i\neq j
\end{cases}$\\
& \\
effective power input & $P_{i}\in{\rm{Unif}}\left([0,10]\right)$\\
& \\
coupling coefficient & $k_{ij}\begin{cases}
=0 & \text{for}\; i=j\\
\in{\rm{Unif}}\left([0.7,1.2]\right) & \text{for}\; i\neq j	
\end{cases}$\\
& \\
\hline
\end{tabular}
\vspace*{0.05in}
\caption{Parameters used for the numerical simulation in Sec. \ref{subsec:NumSimSynthetic}. The indices $i,j\in\{1,\hdots,n\}$ where the number of generators $n=50$ in Sec. \ref{subsec:NumSimSynthetic}.}
\label{table:ParamUncertaintiesSynthetic}
\vspace*{-0.35in}
\end{table}

\subsection{Synthetic Test System}\label{subsec:NumSimSynthetic}
To highlight the scalability of the proposed method, we next consider a power network with $n=50$ generators, and propagate the transient joint state PDFs supported over the 100 dimensional state space $\mathbb{T}^{50}\times\mathbb{R}^{50}$. We take the initial joint state PDF at $t=0$ as
\begin{multline}
\rho_{0}\equiv\rho(t=0,\bm{\theta}(0),\bm{\omega}(0)) = \;{\rm{Unif}}\left(\left([0,2\pi)\;\text{rad}\right)^{n}\right) \\
	\qquad\times {\rm{Unif}}\left(\left([-0.1,0.1]\;\text{rad/s}\right)^{n}\right),
\label{rho0synthetic}	
\end{multline}
and following \cite[Sec. 5]{dorfler2012synchronization}, randomly generate the parameters as in Table \ref{table:ParamUncertaintiesSynthetic}. These parameter ranges are consistent with the same found in \cite{sauerpai1998,andersonBook1977,kundurBook1994}.

\begin{figure}[t]
\centering
\includegraphics[width=0.95\linewidth]{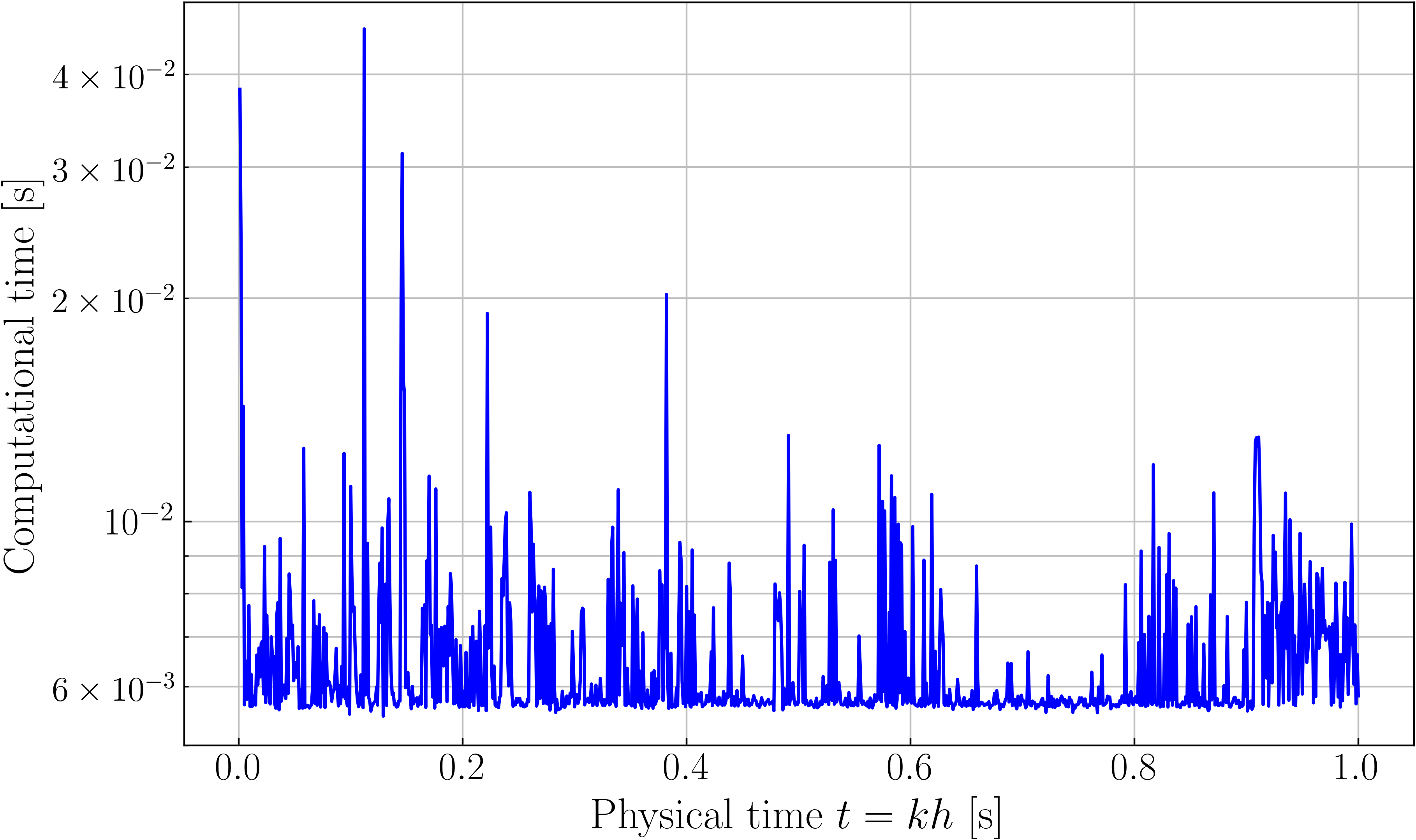}
\caption{\small{The computational time for propagating the transient joint state PDFs over $\mathbb{T}^{50}\times\mathbb{R}^{50}$ for the simulation set up in Sec. \ref{subsec:NumSimSynthetic}.}}
\vspace*{-0.1in}
\label{fig:CompTimeSynthetic}
\end{figure}
\begin{figure}[t]
\centering
\includegraphics[width=0.95\linewidth]{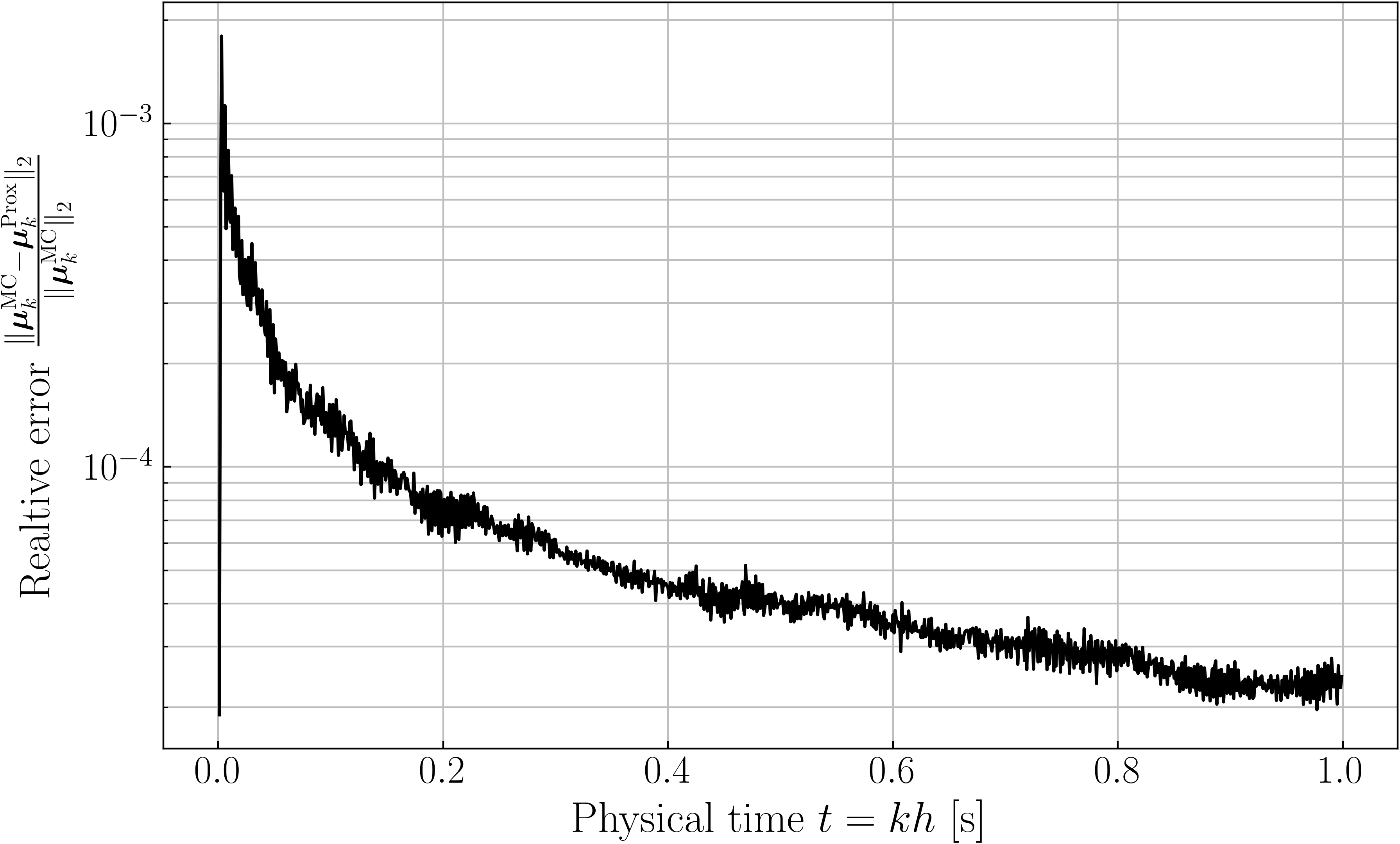}
\caption{\small{The relative error between the empirical mean (i.e., Monte Carlo estimate) $\bm{\mu}_{k}^{\text{MC}}\in\mathbb{T}^{50}\times\mathbb{R}^{50}$, and the mean $\bm{\mu}_{k}^{\text{Prox}}\in\mathbb{T}^{50}\times\mathbb{R}^{50}$ obtained using the proximal updates of the joint state PDFs, for the simulation set up in Sec. \ref{subsec:NumSimSynthetic}.}}
\vspace*{-0.15in}
\label{fig:RelErrSynthetic}
\end{figure}
\begin{figure}[t]
\centering
\includegraphics[width=0.95\linewidth]{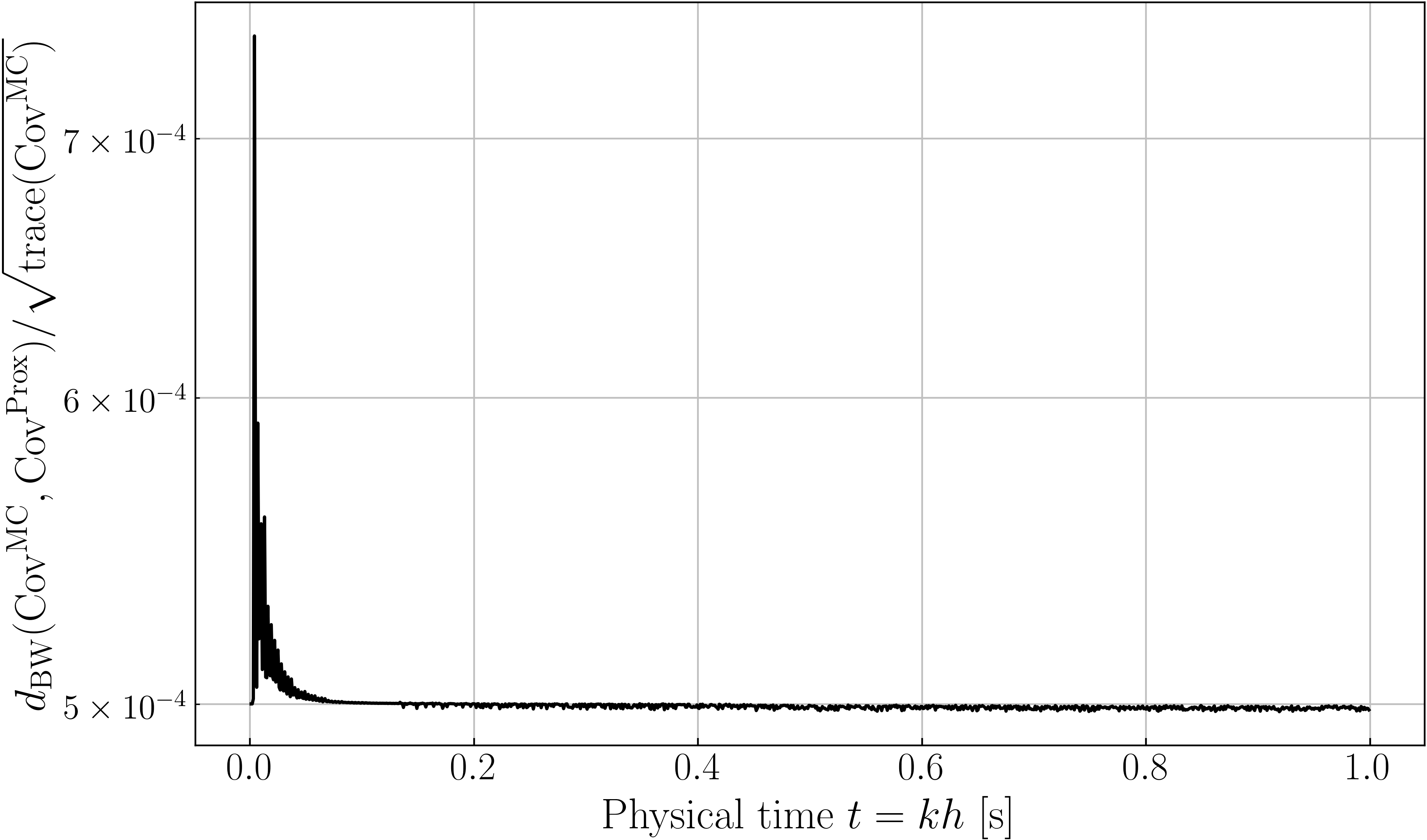}
\caption{\small{The normalized Bures-Wasserstein distance $d_{\rm{BW}}(\cdot,\cdot)$ between the empirical covariance  (i.e., Monte Carlo estimate) ${\rm{Cov}}^{\rm{MC}}\in\mathbb{R}^{100\times 100}$, and the ensemble covariance ${\rm{Cov}}^{\rm{Prox}}\in\mathbb{R}^{100\times 100}$ obtained using the proximal updates of the joint state PDFs, for the simulation set up in Sec. \ref{subsec:NumSimSynthetic}. The normalization is w.r.t. $\sqrt{{\rm{trace}}\left({\rm{Cov}}^{\rm{MC}}\right)}$.}}
\vspace*{-0.2in}
\label{fig:RelErrCovSynthetic}
\end{figure}
\begin{figure}[t]
\centering
\includegraphics[width=0.95\linewidth]{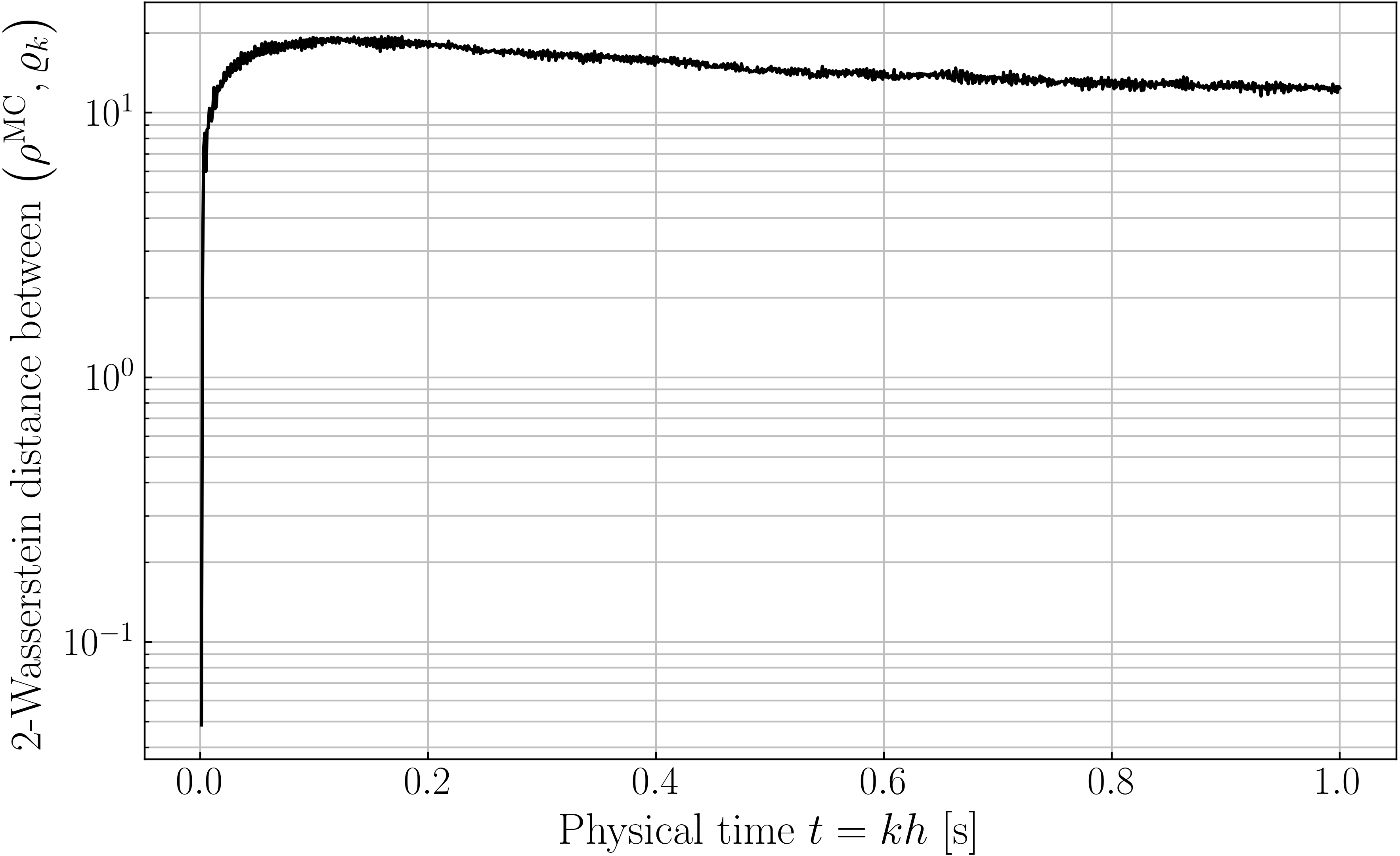}
\caption{\small{The 2-Wasserstein distance between the empirical joint PDF (i.e., Monte Carlo estimate) $\rho^{\rm{MC}}$, and the proximal update $\varrho^{k}$ obtained using the proposed proximal algorithm, for the simulation set up in Sec. \ref{subsec:NumSimSynthetic}.}}
\vspace*{-0.2in}
\label{fig:WassJointPDFSynthetic}
\end{figure}

With $N=2000$ random samples from the initial joint PDF \eqref{rho0synthetic}, and with the aforesaid randomly generated parameters, we employed the procedure detailed in Sec. \ref{sec:ProxAlgo} for propagating $\{\bm{x}_{k}^{i},\varrho_{k}^{i}\}_{i=1}^{N}$. As in Sec. \ref{subsec:NumSimIEEEBus}, we used $h=10^{-3}$, $\varepsilon = 0.05$, $\delta = 10^{-3}$, $\ell_{\max} = 300$ in Algorithm \ref{algo:KineticProx}.

Fig. \ref{fig:CompTimeSynthetic} highlights that for computing high dimensional transient joint state PDFs as in here, the proposed variational framework enjoys remarkably fast computational time. Notice that constructing histograms or other conventional density estimators for the joint state PDF by making a grid over the 100 dimensional state space, is computationally prohibitive in a general purpose computer. Fig. \ref{fig:TotalCompTime} (square markers) shows the total computational times needed for propagating the transient joint state PDFs over the time horizon $[0, 1\;\text{min}]$ for 100 different simulation instances.

To depict the numerical accuracy, Fig. \ref{fig:RelErrSynthetic} plots the (time-varying) relative error between the empirical (i.e., Monte Carlo estimate) mean vector $\bm{\mu}_{k}^{\text{MC}}\in\mathbb{T}^{50}\times\mathbb{R}^{50}$, and the ``proximal mean" vector $\bm{\mu}_{k}^{\text{Prox}}\in\mathbb{T}^{50}\times\mathbb{R}^{50}$. The latter was computed using the proximal updates $\{\bm{x}_{k}^{i},\varrho_{k}^{i}\}_{i=1}^{N}$. Likewise, Fig. \ref{fig:RelErrCovSynthetic} plots the time-varying normalized Bures-Wasserstein distance $d_{\rm{BW}}(\cdot,\cdot)$ between the empirical covariance matrix (i.e., Monte Carlo estimate) ${\rm{Cov}}^{\rm{MC}}\in\mathbb{R}^{100\times 100}$, and the ensemble covariance matrix ${\rm{Cov}}^{\rm{Prox}}\in\mathbb{R}^{100\times 100}$ obtained using the proximal updates. The unnormalized Bures-Wasserstein distance \cite{bhatia2019bures} given by 
\begin{align*}
&d_{\rm{BW}}(\rm{Cov}^{\rm{MC}},\rm{Cov}^{\rm{Prox}}) := \left[\operatorname{trace}\left(\rm{Cov}^{\rm{MC}}+\rm{Cov}^{\rm{Prox}}\right)\right.\\
&\qquad\qquad \left.-2\operatorname{trace}\left(\left(\sqrt{\rm{Cov}^{\rm{MC}}}\:\rm{Cov}^{\rm{Prox}}\:\sqrt{\rm{Cov}^{\rm{MC}}}\right)^{\frac{1}{2}}\right)\right]^{\frac{1}{2}},	
\end{align*}
is a metric on the cone of symmetric positive definite matrices such as the covariance matrices, and can be seen as an absolute error between the covariances. Thus, Fig. \ref{fig:RelErrCovSynthetic} can be viewed as a relative error curve between the covariances. In summary, Figs. \ref{fig:RelErrSynthetic} and \ref{fig:RelErrCovSynthetic} show that the time-varying statistics between the Monte Carlo and the proposed proximal computation remain close at all times.

Since it is not possible to \emph{a priori} determine the dimension of the sufficient statistic of the transient joint state PDFs, Fig. \ref{fig:WassJointPDFSynthetic} shows the 2-Wasserstein metric between the empirical joint PDF (i.e., Monte Carlo estimate) $\rho^{\rm{MC}}$, and the proximal update $\varrho^{k}$. In other words, Fig. \ref{fig:WassJointPDFSynthetic} can be viewed as a nonparametric absolute error curve between $\rho^{\rm{MC}}$ and $\varrho^{k}$. Each snapshot value in Fig. \ref{fig:WassJointPDFSynthetic} is computed by solving the well-known optimal transport linear program (see e.g., \cite[Sec. 4.2-4.3]{halder2014probabilistic}, \cite[Ch. 3]{peyre2019computational}) using the corresponding time-varying weighted point clouds. The dynamics helps guide the concentration of the time-varying point clouds in regions of high and low likelihoods in the state space, and hence as time progresses, the empirical Monte Carlo approximants get closer to the proximal updates. Consequently, Figs. \ref{fig:RelErrSynthetic}, \ref{fig:RelErrCovSynthetic} and \ref{fig:WassJointPDFSynthetic} show that after initial transients, the distances between the Monte Carlo and the proximal statistics decrease with time.


\section{Conclusion}\label{sec:conclusion}
We proposed a nonparametric algorithm to propagate the joint state PDFs subject to networked power system dynamics with stochastic initial conditions, parameters and process noise. Our development is built on recent advances in generalized gradient flows on the space of PDFs, and new analytical results specific to the power system's stochastic dynamics that are presented in this paper. The novelty of the proposed proximal algorithm is that it does not approximate the dynamical or statistical nonlinearities. Instead, the algorithm allows gridless computation by actually exploiting the structure of the PDE governing the evolution of the joint state PDF that is induced by the underlying nonlinear sample path dynamics of a networked power system. Numerical case studies reveal that the computational framework can be a scalable way to propagate stochastic uncertainties in realistic power networks.

Possible future directions include generalizing the proposed method for more complex generator dynamics, and to use this framework to actively steer the uncertainties via optimal control on the space of joint PDFs.

\end{document}